\documentclass{amsart}
\usepackage{amssymb}
\usepackage{amsfonts}

\newtheorem{theorem}{Theorem}
\theoremstyle{plain}

\newtheorem{corollary}[theorem]{Corollary}

\newtheorem{definition}[theorem]{Definition}

\newtheorem{lemma}[theorem]{Lemma}

\newtheorem{proposition}[theorem]{Proposition}
\newtheorem{remark}[theorem]{Remark}

\numberwithin{equation}{section}
\numberwithin{theorem}{section}
\input{tcilatex}

\begin{document}
\title{The Weil representation in characteristic two}
\author{Shamgar Gurevich}
\address{Department of Mathematics, University of California, Berkeley, CA
94720, USA. }
\email{shamgar@math.berkeley.edu}
\author{Ronny Hadani}
\address{Department of Mathematics, University of Chicago, IL, 60637, USA.}
\email{hadani@math.uchicago.edu}
\date{March 24, 2008.}
\thanks{\copyright \ Copyright by S. Gurevich and R. Hadani, March. 24, 2008.
All rights reserved.}

\begin{abstract}
In this paper we construct a new variant of the Weil representation,
associated with a symplectic vector space $\left( V,\omega \right) $ defined
over a finite field of characteristic two. Our variant is a representation $%
\rho :AMp\left( V\right) \rightarrow GL\left( \mathcal{H}\right) $, where
the group $AMp\left( V\right) $ is forth cover of a group $ASp\left(
V\right) $, which is a non-trivial gluing of the symplectic group $Sp\left(
V\right) $ and the dual group $V^{\ast }$. In particular, the group $%
ASp\left( V\right) $ contains Weil's pseudo-symplectic group as a strict
subgroup. In the course, we develop the formalism of canonical vector
spaces, which enables us to realize the group $AMp\left( V\right) $ and the
representation $\rho $ in a transparent manner.
\end{abstract}

\maketitle

\section{\protect \bigskip Introduction}

In his celebrated 1964 Acta paper \cite{W} Andr\`{e} Weil constructed a
distinguished unitary representation $\rho _{\mathrm{Weil}}$, which is
associated with a symplectic vector space $\left( V,\omega \right) $ over a
local field $\mathbb{F}$, now referred to as the \textit{Weil representation}%
. The Weil representation has many fascinating properties which have
gradually been brought to light over the last few decades. It now appears
that this representation is a central object in mathematics: Lying at the
fabric of the theory of harmonic analysis and bridging between various
topics in mathematics and physics, including classical invariant theory, the
theory of theta functions and automorphic forms and last (but probably not
least) quantum mechanics.

In his paper, Andr\`{e} Weil constructs $\rho _{\mathrm{Weil}}$, first, in
the setup when $\mathbb{F}$ is a local field of characteristic zero or a
local field of finite characteristic $p$ where $p\neq 2$. In this case, $%
\rho _{\mathrm{Weil}}$ is a representation of a double cover of the
symplectic group $Sp\left( V\right) $ (called the metaplectic cover). Then,
he proceeds \cite{BL, W} to construct $\rho _{\mathrm{Weil}}$ in the more
intricate setup when $\mathbb{F}$ is a local field of characteristic $2$. In
this case, $\rho _{\mathrm{Weil}}$ is a representation of a double cover of
the \textit{pseudo-symplectic} group $Ps\left( V\right) $, which is a
non-trivial gluing of an orthogonal group and the dual vector space $V^{\ast
}$, i.e., we have a short exact sequence of groups 
\begin{equation*}
1\rightarrow V^{\ast }\rightarrow Ps\left( V\right) \rightarrow O\left(
Q\right) \rightarrow 1,
\end{equation*}%
where, $Q:V\rightarrow \mathbb{\mathbb{F}}$ is the quadratic form $Q\left(
v\right) =\beta \left( v,v\right) $ for some non-symmetric bilinear form $%
\beta :V\times V\rightarrow $ $\mathbb{\mathbb{F}}$ such that $\beta \left(
v,u\right) -\beta \left( u,v\right) =\omega \left( v,u\right) $, and $%
O\left( Q\right) $ is the associated orthogonal group.

A comparison between the constructions in the two setups suggests that the
definition of $\rho _{\mathrm{Weil}}$ when $\mathbb{F}$ is a field of
characteristic $2$ is unsatisfactory for the reason that the
pseudo-symplectic group is not related anymore to the symplectic group.
Instead, $Ps\left( V\right) $ surjects onto the orthogonal subgroup\footnote{%
Accidently, in characteristic two, the orthogonal group $O\left( Q\right) $
appears as a subgroup of the symplectic group $Sp\left( V\right) $.} $%
O\left( Q\right) \subsetneqq Sp\left( V\right) $ which, as an algebraic
group, is of smaller dimension 
\begin{eqnarray*}
\dim Sp\left( V\right)  &=&\dim V\left( \dim V+1\right) /2, \\
\dim O\left( Q\right)  &=&\dim V\left( \dim V-1\right) /2\text{.}
\end{eqnarray*}

A natural question to pose at this point is whether there exists an
extension of $\rho _{\mathrm{Weil}}$ to a representation $\rho $ (which acts
on the same Hilbert space) of a larger group $G$, which contains the
pseudo-symplectic group as a subgroup and surjects onto the symplectic group.

\subsection{Main results}

In this paper we construct a new variant of the Weil representation
associated to a symplectic vector space $\left( V,\omega \right) $ defined
over a finite field $\mathbb{F}$ of characteristic two.

\subsubsection{Projective Weil representation}

We describe a group $ASp\left( V\right) $, that we call the affine
symplectic group, which contains the pseudo-symplectic group $Ps\left(
V\right) $ as a subgroup and is glued from the symplectic group $Sp\left(
V\right) $ and the dual abelian group $V^{\ast }$ 
\begin{equation*}
1\rightarrow V^{\ast }\rightarrow ASp\left( V\right) \rightarrow Sp\left(
V\right) \rightarrow 1.
\end{equation*}

In addition, we construct a projective Weil representation 
\begin{equation*}
\widetilde{\rho }:ASp\left( V\right) \rightarrow PGL\left( \mathcal{H}%
\right) ,
\end{equation*}
which extends as a projective representation the representation $\rho _{%
\mathrm{Weil}}$.

\subsubsection{Linear Weil representation}

We describe a group $AMp\left( V\right) ,$ that we call the affine
metaplectic group, which is a central extension of $ASp\left( V\right) $ by
the group $\mu _{4}$ of 4th roots of unity. In addition, we construct a
linear Weil representation 
\begin{equation*}
\rho :AMp\left( V\right) \rightarrow GL\left( \mathcal{H}\right) ,
\end{equation*}
lying over the projective representation $\widetilde{\rho }$.

\subsubsection{Splitting of the Weil representation}

We describe a splitting homomorphism $s:Mp(\widetilde{V})\rightarrow
AMp\left( V\right) $ and consequently a pull-back representation 
\begin{equation*}
\rho _{\widetilde{V}}=\rho \circ s:Mp(\widetilde{V})\rightarrow GL\left( 
\mathcal{H}\right) \text{,}
\end{equation*}%
where $(\widetilde{V},\widetilde{\omega })$ is a free symplectic module over
the the ring $W_{2}\left( \mathbb{F}\right) $ of (level 2) truncated Witt
vectors which reduces to $\left( V,\omega \right) $ mod $2$ and $Mp(%
\widetilde{V})$ is a central extension of the symplectic group $Sp(%
\widetilde{V})$ by the group $\mu _{2}=\left \{ \pm 1\right \} $.

\subsubsection{The formalism of canonical vector spaces}

In the course, we develop the formalism of canonical vector spaces, which
enables us to realize the group $AMp\left( V\right) $, the representations $%
\rho $ and the splitting homomorphism $s:Mp(\widetilde{V})\rightarrow
AMp\left( V\right) $ in a transparent manner and, moreover, serves as an
appropriate conceptual framework for the study of these objects. The
development of this formalism in the characteristic two setting constitutes
the main technical contribution of this paper.

We devote the rest of the introduction to a more detailed account of the
main constructions and results of this paper. For simplicity let us assume
that $\mathbb{F}=\mathbb{F}_{2}$.

\subsection{The Heisenberg group\label{Heisenberg_sub}}

Considering $V$ as an abelian group and a bi-additive form $\beta :V\times
V\rightarrow $ $%
\mathbb{Z}
/4%
\mathbb{Z}
$ such that $\beta \left( u,v\right) -\beta \left( v,u\right) =2\omega
\left( u,v\right) \in 
\mathbb{Z}
/4%
\mathbb{Z}
$, one can associate to $\left( V,\beta \right) $ a central extension $%
H\left( V\right) =H_{\beta }\left( V\right) $ 
\begin{equation*}
0\rightarrow 
\mathbb{Z}
/4%
\mathbb{Z}
\rightarrow H\left( V\right) \rightarrow V\rightarrow 0\text{,}
\end{equation*}%
called the \textit{Heisenberg group}. The group of automorphisms of $H\left(
V\right) $ acting trivially on the center, denoted by $ASp\left( V\right) $
and refer to as the \textit{affine symplectic group}, fits into a non-split
exact sequence%
\begin{equation*}
0\rightarrow V^{\ast }\rightarrow ASp\left( V\right) \rightarrow Sp\left(
V\right) \rightarrow 1.
\end{equation*}

Concretely, the elements of $ASp\left( V\right) $ can be presented as pairs $%
\left( g,\alpha \right) $ where $g\in Sp\left( V\right) $ and $\alpha
:V\rightarrow 
\mathbb{Z}
/4%
\mathbb{Z}
$ satisfies the condition 
\begin{equation}
\alpha \left( v_{1}+v_{2}\right) -\alpha \left( v_{1}\right) -\alpha \left(
v_{2}\right) =\beta \left( gv_{1},gv_{2}\right) -\beta \left(
v_{1},v_{2}\right) \text{.}  \label{condition1_eq}
\end{equation}

\subsubsection{Weil's Heisenberg group}

The previous development should be contrasted with the standard construction
that appears in \cite{W}. There, the Heisenberg group is an extension of $V$
by the field $\mathbb{F}$, which is associated to a bilinear form $\beta
:V\times V\rightarrow $ $\mathbb{F}$ such that $\beta \left( u,v\right)
-\beta \left( v,u\right) =\omega \left( u,v\right) $. The group of
automorphisms which act trivially on the center is the pseudo-symplectic
group $Ps\left( V\right) $ whose elements can be presented as pairs $\left(
g,\alpha \right) $ where $g\in Sp\left( V\right) $ and $\alpha :V\rightarrow 
\mathbb{F}$ satisfies the following (polarization) condition: 
\begin{equation}
\alpha \left( v_{1}+v_{2}\right) -\alpha \left( v_{1}\right) -\alpha \left(
v_{2}\right) =\beta \left( gv_{1},gv_{2}\right) -\beta \left(
v_{1},v_{2}\right) \text{.}  \label{condition2_eq}
\end{equation}

The reason why the pseudo-symplectic group $Ps\left( V\right) $ is strictly
smaller than the affine symplectic group $ASp\left( V\right) $ is because
Equation (\ref{condition1_eq}) admits solutions\ for every $g\in Sp\left(
V\right) $, while Equation (\ref{condition2_eq}) admits solutions only when $%
g$ lies in the orthogonal group $O\left( Q\right) \subset Sp\left( V\right) $
where $Q\left( v\right) =\beta \left( v,v\right) $.

This phenomena can be appreciated already in the following simplified
situation: Let $\beta :$ $\mathbb{F\times }$ $\mathbb{F\rightarrow }$ $%
\mathbb{F}$ be the bilinear form given by 
\begin{equation*}
\beta \left( x,y\right) =xy\text{.}
\end{equation*}

It can be easily verified that there is no quadratic form $\alpha :$ $%
\mathbb{F\rightarrow }$ $\mathbb{F}$ which polarizes $\beta $, namely,
satisfies the condition 
\begin{equation*}
\alpha \left( x+y\right) -\alpha \left( x\right) -\alpha \left( y\right)
=\beta \left( x,y\right) .
\end{equation*}

However, considering the quadratic form $\widetilde{\alpha }:%
\mathbb{Z}
/4%
\mathbb{Z}
\rightarrow 
\mathbb{Z}
/4%
\mathbb{Z}
$ given by $\widetilde{\alpha }\left( x\right) =x^{2}$, it can be easily
verified that $\widetilde{\alpha }$ descends to a function $\alpha :\mathbb{%
F\rightarrow }%
\mathbb{Z}
/4%
\mathbb{Z}
$ which satisfies 
\begin{equation*}
\alpha \left( x+y\right) -\alpha \left( x\right) -\alpha \left( y\right)
=2\beta \left( x,y\right) \in 
\mathbb{Z}
/4%
\mathbb{Z}
\text{.}
\end{equation*}

\subsection{The Heisenberg representation}

Fixing a faithful character $\psi $ of the center of $H\left( V\right) $,
there exists a unique irreducible representation $\pi :H\left( V\right)
\rightarrow GL\left( \mathcal{H}\right) $ with a central character $\psi $ -
this is the \textit{Stone-von Neumann property} (S-vN property for short).
We refer to this representation as the \textit{Heisenberg representation}

\subsection{Realizations of the Heisenberg representation}

The Heisenberg representation admits a special family of models
(realizations) associated with \textit{enhanced Lagrangian} subspaces in $V$.

An enhanced Lagrangian is a homomorphism section (with respect to the
natural projection $H\left( V\right) \rightarrow V$) $\tau :L\rightarrow
H\left( V\right) $, where $L\in Lag\left( V\right) $ is a Lagrangian
subspace in $V$. For such a $\tau $ one can associate a model $\left( \pi
_{L},H\left( V\right) ,\mathcal{H}_{L}\right) $ (abusing the notations a
bit) of the Heisenberg representation, which is defined as follows:

The vector space $\mathcal{H}_{L}$ consists of functions $f:H\left( V\right)
\rightarrow 
\mathbb{C}
$ satisfying 
\begin{equation*}
f\left( z\cdot \tau \left( l\right) \cdot h\right) =\psi \left( z\right)
f\left( h\right) ,
\end{equation*}
with $l\in L$ and $z\in Z\left( H\left( V\right) \right) $ a central element
and the action $\pi _{L}$ is given by right translations.

The collection of models $\left \{ \mathcal{H}_{L}\right \} $ forms a vector
bundle $\mathfrak{H}$ on the set $ELag\left( V\right) $ of enhanced
Lagrangians, with fibers $\mathfrak{H}_{L}=\mathcal{H}_{L}$.

\subsection{The strong Stone-von Neumann property}

A basic technical statement is a strong variant of the Stone-von Neumann
property and it asserts that the vector bundle $\mathfrak{H}^{\otimes 4}$
admits a \textbf{natural} trivialization, that is, existence of a canonical
system of intertwining morphisms $T_{M,L}:\mathcal{H}_{L}^{\otimes
4}\rightarrow \mathcal{H}_{M}^{\otimes 4}$ which satisfies the following
multiplicativity condition:%
\begin{equation*}
T_{N,M}\circ T_{M,L}=T_{N,L}\text{,}
\end{equation*}%
for every $N,M,L\in ELag\left( V\right) $.

\subsection{The Weil gerbe}

There exists a groupoid category $\mathcal{W}$ which is naturally associated
with the vector bundle $\mathfrak{H}$ and encoding in its structure the
strong S-vN property:

\begin{itemize}
\item An object in $\mathcal{W}$ is a triple $\left( \mathfrak{E,}%
\{E_{M,L}\},\varphi \right) $, where $\mathfrak{E}$ is a vector bundle on $%
ELag\left( V\right) $, $\{E_{M,L}\}$ is a trivialization of $\mathfrak{E}$
and $\varphi :\mathfrak{E}\overset{\simeq }{\rightarrow }\mathfrak{H}$ an
isomorphism of vector bundles which satisfies that $\varphi ^{\otimes 4}:%
\mathfrak{E}^{\otimes 4}\overset{\simeq }{\rightarrow }\mathfrak{H}^{\otimes
4}$ is an isomorphism of trivialized vector bundles.

\item A morphism in $\mathcal{W}$ is an isomorphism of trivialized vector
bundles $f:\mathfrak{E}_{1}\overset{\simeq }{\rightarrow }\mathfrak{E}_{2}$
which satisfies $\varphi _{2}^{\otimes 4}\circ f^{\otimes 4}=\varphi
_{1}^{\otimes 4}$.
\end{itemize}

\textbf{Main property: }The groupoid $\mathcal{W}$ is a gerbe with band $\mu
_{4}$ which means that every two objects in $\mathcal{W}$ are isomorphic and 
\textrm{Mor}$_{\mathcal{W}}\left( \mathfrak{E,E}\right) \simeq \mu _{4}$,
for every $\mathfrak{E\in }\mathcal{W}$.

We refer to $\mathcal{W}$ as the \textit{Weil \ }gerbe.

\subsubsection{Action of the affine symplectic group}

The group $ASp\left( V\right) $ naturally acts on the Weil gerbe. The action
of each element $g\in ASp\left( V\right) $ \ is given by the pull-back
functor $g^{\ast }:\mathcal{W\rightarrow W}$, sending a vector bundle $%
\mathfrak{E}$ to its pull-back $g^{\ast }\mathfrak{E}$.

\subsection{Canonical vector space}

There exists a tautological fiber functor $\Gamma $ from $\mathcal{W}$ to
the category $\mathsf{Vect}$ of complex vector spaces, sending an object $%
\mathfrak{E\in }\mathcal{W}$ to the space of "horizontal sections" $\Gamma
_{hor}\left( ELag\left( V\right) ,\mathfrak{E}\right) $, consisting of
systems 
\begin{equation*}
\left( f_{L}\in \mathfrak{E}_{L}:L\in ELag\left( V\right) \right) ,
\end{equation*}%
such that $E_{M,L}\left( f_{L}\right) =f_{M}$, for every $M,L\in ELag\left(
V\right) $.

In addition, by general considerations \cite{B}, there is a central
extension 
\begin{equation*}
1\rightarrow \mu _{4}\rightarrow AMp\left( V\right) \rightarrow ASp\left(
V\right) \rightarrow 1,
\end{equation*}%
which is naturally associated with the action of $ASp\left( V\right) $ on
the groupoid $\mathcal{W}$: An element of $AMp\left( V\right) $ is a pair $%
\left( g,\iota \right) $ where $g\in ASp\left( V\right) $ and $\iota
:g^{\ast }\overset{\simeq }{\rightarrow }Id$ is an isomorphism of functors.

Finally, there is a natural homomorphism 
\begin{equation*}
\rho :AMp\left( V\right) \rightarrow \mathrm{Aut}\left( \Gamma \right) .
\end{equation*}

\textbf{Summary}: The fundamental object is the Weil gerbe $\mathcal{W}$
equipped with the action of the affine symplectic group $ASp\left( V\right) $%
. The following structures are canonically associated with this object:

\begin{itemize}
\item A group $AMp\left( V\right) $, which is a central extension of $%
ASp\left( V\right) $ by the group $\mu _{4}$.

\item A fiber functor $\Gamma :\mathcal{W}$ $\rightarrow \mathsf{Vect}$,
which might be thought of as a vector space twisted by the gerbe $\mathcal{W}
$.

\item A representation $\rho $ $:AMp\left( V\right) $ $\rightarrow $ $%
\mathrm{Aut}\left( \Gamma \right) $, which we refer to as the canonical
model of the Weil representation.
\end{itemize}

After choosing an object $\mathfrak{E}$ $\in \mathcal{W}$ one returns to a
more traditional setting, obtaining a homomorphism 
\begin{equation*}
\rho _{\mathfrak{E}}:AMp\left( V\right) \rightarrow GL\left( \Gamma \left( 
\mathfrak{E}\right) \right) .
\end{equation*}

\subsection{Splitting of the Weil representation}

The splitting depends on an auxiliary data of a free symplectic module $(%
\widetilde{V},\widetilde{\omega })$ over the ring $W_{2}\left( \mathbb{F}%
\right) $ which reduces to $\left( V,\omega \right) $ mod $2$.

The following structures are naturally associated with such data:

\begin{itemize}
\item A homomorphism $Sp(\widetilde{V})\rightarrow ASp\left( V\right) $.

\item A gerbe $\mathcal{W}^{s}$ with band $\mu _{2}$, equipped with an
action of the group $Sp(\widetilde{V})$.

\item A faithful (splitting) functor $S:\mathcal{W}^{s}\rightarrow \mathcal{W%
}$ which is compatible with the action of the group $Sp(\widetilde{V})$ on
both sides.
\end{itemize}

In complete analogy to the definition of the group $AMp\left( V\right) $,
there is a central extension 
\begin{equation*}
1\rightarrow \mu _{2}\rightarrow Mp(\widetilde{V})\rightarrow Sp(\widetilde{V%
})\rightarrow 1,
\end{equation*}%
which is naturally associated with the action of $Sp(\widetilde{V})$ on the
groupoid $\mathcal{W}^{s}$.

In addition, the splitting functor $S:\mathcal{W}^{s}\rightarrow \mathcal{W}$
yields a homomorphism 
\begin{equation*}
s:Mp(\widetilde{V})\rightarrow AMp\left( V\right) ,
\end{equation*}%
and consequently a representation $\rho _{\widetilde{V}}:Mp(\widetilde{V}%
)\rightarrow \mathrm{Aut}\left( \Gamma \right) .$

\subsection{Structure of the paper}

Apart from the introduction, the paper consists of five sections and an
appendix.

\begin{itemize}
\item In Section \ref{Weil_sec}, we introduce the Weil representation
associated with a symplectic vector space $\left( V,\omega \right) $ defined
over a field of characteristic $2$. We begin by describing an appropriate
Heisenberg group $H\left( V\right) $. Then we describe the group $ASp\left(
V\right) $ of automorphisms of $H\left( V\right) $ which act trivially on
the center. We proceed to describe the Heisenberg representation and
formulate the Stone-von Neumann property for this representation. We end
this section with two theorems. The first Theorem (Theorem \ref{Weilrep_thm}%
) asserts the existence of the Weil representation $\rho $ of the affine
metaplectic group $AMp\left( V\right) $. The second Theorem (Theorem \ref%
{Weilrep-split_thm}) asserts the existence of a splitting of $\rho $ over
the group $Sp(\widetilde{V})$ which amounts to a representation of the
metaplectic group $Mp(\widetilde{V})$.

\item In Section \ref{canonical_sec}, we develop the formalism of canonical
vector spaces. We begin by describing a special family of models of the
Heisenberg representation which are associated with enhanced Lagrangian
subspaces in $V$ and the associated Heisenberg vector bundle $\mathfrak{H}$
on the set $ELag\left( V\right) $ of enhanced Lagrangians. We proceed to
define the notion of a trivialization of an Heisenberg vector bundle. The
main statement is Theorem \ref{SS-vN_thm} asserting the existence of a
canonical trivialization of the vector bundle $\mathfrak{H}^{\otimes 4}$.
Using Theorem \ref{SS-vN_thm}, we define the Weil gerbe $\mathcal{W}$, which
is then used to construct the canonical model of the Weil representation,
proving, in particular, Theorem \ref{Weilrep_thm}$.$

\item In section \ref{S-vN_sec}, we describe the construction of the
canonical trivialization asserted in Theorem \ref{SS-vN_thm}, specifically,
we describe the canonical intertwining morphisms between transversal models
of the Heisenberg representation and the cocycle $C$ which is associated
with them. The main Theorem is Theorem \ref{power_thm} which asserts that $%
C^{4}=1$. The rest of the section is devoted to the proof of Theorem \ref%
{power_thm}. In the course, we obtain some results in the theory of
symmetric spaces over the ring $%
\mathbb{Z}
/4%
\mathbb{Z}
$.

\item In Section \ref{canonical-split_sec}, the formalism of canonical
vector spaces is further developed. We begin by introducing the notion of an
oriented Lagrangian in $\widetilde{V}$. We then describe an Heisenberg
vector bundle $\widetilde{\mathfrak{H}}$ on the set $OLag(\widetilde{V})$ of
oriented Lagrangians in $\widetilde{V}$. The main statement is Theorem \ref%
{S-vN-split_thm} asserting the existence of a natural trivialization of the
vector bundle $\widetilde{\mathfrak{H}}^{\otimes 2}$. Using Theorem \ref%
{S-vN-split_thm}, we define the splitting of the Weil gerbe $S:\mathcal{W}%
^{s}\rightarrow \mathcal{W}$, which is then used to prove Theorem \ref%
{Weilrep-split_thm}$.$

\item In section \ref{S-vN-split_sec}, we describe the construction of the
canonical trivialization asserted in Theorem \ref{S-vN-split_thm},
specifically, we describe a natural normalization of the canonical
intertwining morphisms and the cocycle $C$ which is associated with these
normalized intertwining morphisms. The main Theorem of this section is
Theorem \ref{cocycle_thm} which asserts that $C^{2}=1$.

\item In Appendix \ref{proofs_sec}, we give the proofs of all technical
statements which appear in the body of the paper.
\end{itemize}

\subsection{Acknowledgements}

We would like to thank our teacher J. Bernstein for his interest and
guidance. It is a pleasure to thank D. Kazhdan for sharing with us his
thoughts about the nature of canonical Hilbert spaces. We are also thankful
to P. Etingof for telling us about the peculiarities of the Weil
representation in even characteristic, a remark which initiated this
direction of research. We thank M. Nori for explaining to us several useful
ideas from algebraic number theory. Finally, we thank T. Schedler for many
stimulating discussions about Heisenberg groups, Weil representations and
commutative rings with even number of elements.

\section{The Weil representation \label{Weil_sec}}

\subsection{General setting}

\subsubsection{Fields and rings}

Let $K$ be an unramified extension of degree $d$ of the the 2-adic
completion $%
\mathbb{Q}
_{2}$. Let $\mathcal{O}_{K}\subset K$ be the ring of integers, $\mathfrak{m}%
_{K}\subset \mathcal{O}_{K}$ the unique maximal ideal with its standard
generator $2\in \mathfrak{m}_{K}$ and $k=\mathcal{O}_{K}/\mathfrak{m}_{K}$
the residue field, $k=\mathbb{F}_{2^{d}}$. Finally we denote by $R$ the ring 
$\mathcal{O}_{K}/\mathfrak{m}_{K}^{2}$ and remind that we have the trace map%
\begin{equation*}
tr:R\rightarrow 
\mathbb{Z}
/4\text{.}
\end{equation*}

\subsubsection{Symplectic module}

Let $(\widetilde{V},\widetilde{\omega })$ be a free symplectic module over $%
R $ of rank $2n$. Let $V=\widetilde{V}/\mathfrak{m}_{K}$ be the quotient $k$%
-vector space. The form $\omega =2\widetilde{\omega }$ factors to give a
non-degenerate skew symmetric form on $V$ with values in $R$. We denote by $%
Sp(\widetilde{V})=Sp(\widetilde{V},\widetilde{\omega })$ and by $Sp\left(
V\right) =Sp\left( V,\omega \right) $ the corresponding groups of linear
symplectomorphisms.

\subsubsection{The cocycle associated with a Lagrangian splitting}

Let $\widetilde{S}$ be a Lagrangian splitting $\widetilde{V}=\widetilde{L}%
\times \widetilde{M}$; we define a bilinear form $\widetilde{\beta }=%
\widetilde{\beta }_{\widetilde{S}}:\widetilde{V}\times \widetilde{V}%
\rightarrow R,$ given by $\widetilde{\beta }((\widetilde{l}_{1},\widetilde{m}%
_{1}),(\widetilde{l}_{2},\widetilde{m}_{2}))=\widetilde{\omega }(\widetilde{l%
}_{1},\widetilde{m}_{2}),$ for $\widetilde{l}_{i}\in \widetilde{L}$ and $%
\widetilde{m}_{i}\in \widetilde{M}$, $i=1,2$. A direct verification reveals
that $\widetilde{\beta }$ is a cocycle, namely 
\begin{equation*}
0=d\widetilde{\beta }(\widetilde{v}_{1},\widetilde{v}_{2},\widetilde{v}_{3})=%
\widetilde{\beta }(\widetilde{v}_{2},\widetilde{v}_{3})-\widetilde{\beta }(%
\widetilde{v}_{1}+\widetilde{v}_{2},\widetilde{v}_{3})+\widetilde{\beta }(%
\widetilde{v}_{1},\widetilde{v}_{2}+\widetilde{v}_{3})-\widetilde{\beta }(%
\widetilde{v}_{1},\widetilde{v}_{2})\text{,}
\end{equation*}%
for every $\widetilde{v}_{1},\widetilde{v}_{2},\widetilde{v}_{3}\in 
\widetilde{V}$. In addition, $\widetilde{\beta }(\widetilde{v}_{1},%
\widetilde{v}_{2})-\widetilde{\beta }(\widetilde{v}_{2},\widetilde{v}_{1})=%
\widetilde{\omega }\left( \widetilde{v}_{1},\widetilde{v}_{2}\right) $, for
every $\widetilde{v}_{1},\widetilde{v}_{2}\in \widetilde{V}$. Finally, we
consider the form $\beta =2\widetilde{\beta }$, which factors to give a
cocycle on $V$ with values in $R$, with the property that $\beta \left(
v_{1},v_{2}\right) -\beta \left( v_{2},v_{1}\right) =\omega \left(
v_{1},v_{2}\right) $.

For the rest of this paper we fix a splitting $\widetilde{S}$ and denote by $%
\widetilde{\beta }$ and $\beta $ the corresponding cocycles.

\subsection{The Heisenberg group}

Considering $V$ as an abelian group, we can associate to the pair $\left(
V,\beta \right) $ a central extension 
\begin{equation*}
0\rightarrow R\rightarrow H_{\beta }\left( V\right) \rightarrow V\rightarrow
0\text{.}
\end{equation*}

The group $H\left( V\right) =H_{\beta }\left( V\right) $ is called the 
\textit{Heisenberg group} associated with the cocycle $\beta $. More
concretely, the Heisenberg group can be presented as $H\left( V\right)
=V\times R$, with the multiplication given by 
\begin{equation*}
\left( v_{1},z_{1}\right) \cdot \left( v_{2},z_{2}\right) =\left(
v_{1}+v_{2},z_{1}+z_{2}+\beta \left( v_{1},v_{2}\right) \right) \text{. }
\end{equation*}

The center of $H\left( V\right) $ is $Z=Z_{H\left( V\right) }=\left \{
\left( 0,z\right) :z\in R\right \} $.

\subsubsection{Automorphisms of the Heisenberg group}

Let us denote by $ASp\left( V\right) $ the group of automorphisms of $%
H\left( V\right) $ acting trivially on the center. The group $ASp\left(
V\right) $ can be presented as follows: Given an element $g\in Sp\left(
V\right) $, we denote by $\Sigma _{g}$ the set consisting of "quadratic
functions" $\alpha :V\rightarrow R$, which satisfies 
\begin{equation*}
\alpha \left( v_{1}+v_{2}\right) -\alpha \left( v_{1}\right) -\alpha \left(
v_{2}\right) =\beta \left( g\left( v_{1}\right) ,g\left( v_{2}\right)
\right) -\beta \left( v_{1},v_{2}\right) \text{,}
\end{equation*}%
for every $v_{1},v_{2}\in V$. We can write $ASp\left( V\right) =\left \{
\left( g,\alpha \right) :g\in Sp\left( V\right) ,\alpha \in \Sigma
_{g}\right \} $, with the multiplication rule given by 
\begin{equation*}
\left( g,\alpha _{g}\right) \cdot \left( h,\alpha _{h}\right) =\left( g\cdot
h,Ad_{h^{-1}}\left( \alpha _{g}\right) +\alpha _{h}\right) \text{,}
\end{equation*}%
where $Ad_{h^{-1}}\left( \alpha _{g}\right) \left( v\right) =\alpha
_{g}\left( h\left( v\right) \right) $, for every $v\in V$. \ An element $%
\left( g,\alpha \right) \in ASp\left( V\right) $ \ defines the automorphism $%
\left( v,z\right) \longmapsto \left( g\left( v\right) ,z+\alpha \left(
v\right) \right) $ of $H\left( V\right) .$

The group $ASp\left( V\right) $ fits into a non-split exact sequence%
\begin{equation*}
0\rightarrow V^{\vee }\rightarrow ASp\left( V\right) \rightarrow Sp\left(
V\right) \rightarrow 1,
\end{equation*}%
where $V^{\vee }$ is the dual group $V^{\vee }=\mathrm{Hom}\left( V,R\right) 
$. We will refer to $ASp\left( V\right) $ as the \textit{affine symplectic
group}.

It is instructive to give an explicit description of an element in $\Sigma
_{g}$. Let $\widetilde{g}\in Sp(\widetilde{V})$ be an element which lies
over $g$ and let $\widetilde{\alpha }:\widetilde{V}\rightarrow R$ be the
quadratic form defined by $\widetilde{\alpha }\left( \widetilde{v}\right) =%
\widetilde{\beta }\left( \widetilde{g}\left( \widetilde{v}\right) ,%
\widetilde{g}\left( \widetilde{v}\right) \right) -\widetilde{\beta }\left( 
\widetilde{v},\widetilde{v}\right) $.

\begin{lemma}
\label{tech1_lemma}The quadratic form $\widetilde{\alpha }$ factors to a
function $\alpha _{\widetilde{g}}:V\rightarrow R$, moreover $\alpha _{%
\widetilde{g}}\in \Sigma _{g}$.
\end{lemma}

For a proof, see Appendix \ref{proofs_sec}.

\subsection{The Heisenberg representation}

One of the most important attributes of the group $H\left( V\right) $ is
that it admits, principally, a unique irreducible representation - this is
the Stone-von Neumann property (S-vN property for short). The precise
statement goes as follows. Let $\psi :R\rightarrow 
\mathbb{C}
^{\times }$ be a faithful character.

\begin{theorem}[Stone-von Neumann property]
\label{S-vN_thm}There exists a unique (up to a non-unique isomorphism)
irreducible representation $\left( \pi ,H\left( V\right) ,\mathcal{H}\right) 
$, with central character $\psi $, i.e., $\pi _{|Z}=\psi \cdot Id_{\mathcal{H%
}}$.
\end{theorem}

The representation $\pi $ which appears in the above theorem will be called
the \textit{Heisenberg representation} associated with the central character 
$\psi $. For the rest of this paper we take $\psi \left( z\right) =e^{\frac{%
2\pi i}{4}tr\left( z\right) }$.

\subsection{The Weil representation}

A direct consequence of Theorem \ref{S-vN_thm} is the existence of a
projective representation $\widetilde{\rho }:ASp\left( V\right) \rightarrow
PGL(\mathcal{H)}$. The construction of $\widetilde{\rho }$ out of the
Heisenberg representation $\pi $ is rather standard and it goes as follows.
Considering the Heisenberg representation $\pi $ and an element $g\in
ASp\left( V\right) $, one can define a new representation $\pi ^{g}$ acting
on the same Hilbert space via $\pi ^{g}\left( h\right) =\pi \left( g\left(
h\right) \right) $. Clearly both $\pi $ and $\pi ^{g}$ have the same central
character $\psi $ hence by Theorem \ref{S-vN_thm} they are isomorphic. Since
the space $\mathrm{Hom}_{H\left( V\right) }(\pi ,\pi ^{g})$ is
one-dimensional, choosing for every $g\in ASp\left( V\right) $ a non-zero
representative $\widetilde{\rho }(g)\in \mathrm{Hom}_{H\left( V\right) }(\pi
,\pi ^{g})$ gives the required projective representation. In more concrete
terms, the projective representation $\widetilde{\rho }$ is characterized by
the formula 
\begin{equation}
\widetilde{\rho }\left( g\right) \pi \left( h\right) \widetilde{\rho }\left(
g^{-1}\right) =\pi \left( g\left( h\right) \right) ,  \label{Egorov}
\end{equation}%
for every $g\in ASp\left( V\right) $ and $h\in H\left( V\right) $.

Our goal is to prove the following theorem

\begin{theorem}[The Weil representation]
\label{Weilrep_thm}There exists a group $AMp\left( V\right) $ which is a
central extension of $ASp\left( V\right) $ by the group $\mu _{4}$ of 4th
roots of unity and a linear representation $\rho :AMp\left( V\right)
\rightarrow GL\left( \mathcal{H}\right) $ lying over $\widetilde{\rho }$.
\end{theorem}

For a proof, see Section \ref{canonical_sec}.

\subsection{Splitting of the Weil representation}

There is a natural homomorphism $Sp(\widetilde{V})\rightarrow ASp\left(
V\right) $, sending an element $\widetilde{g}\in Sp(\widetilde{V})$ to the
element $\left( g,\alpha _{\widetilde{g}}\right) \in ASp\left( V\right) $
where $g\in Sp\left( V\right) $ is the reduction of $\widetilde{g}$ mod 2
and $\alpha _{\widetilde{g}}\in \Sigma _{g}$ is the "quadratic function"
associated to the lift $\widetilde{g}\mapsto g$ (see Theorem \ref%
{tech1_lemma}). Our goal is to prove the following splitting theorem:

\begin{theorem}
\label{Weilrep-split_thm}There exists a group $Mp(\widetilde{V})$ which is a
central extension of $Sp(\widetilde{V})$ by the group $\mu _{2}=\left \{ \pm
1\right \} $ and a homomorphism of central extensions $s:$ $Mp(\widetilde{V}%
)\rightarrow AMp\left( V\right) $.
\end{theorem}

For a proof, see Section \ref{canonical-split_sec}.

As a direct consequence of Theorem \ref{Weilrep-split_thm} we obtain the
pull-back representation 
\begin{equation*}
\rho _{\widetilde{V}}=\rho \circ s:Mp(\widetilde{V})\rightarrow GL\left( 
\mathcal{H}\right) .
\end{equation*}

\section{Canonical vector space\label{canonical_sec}}

\subsection{Models of the Heisenberg representation}

Although, the representation $\pi $ is unique, it admits a multitude of
different models (realizations), in fact, this is one of its most
interesting and powerful attributes. In this work we will be interested in a
particular family of such models, which are associated with \textit{enhanced
Lagrangian} subspaces in $V$.

\subsubsection{Enhanced Lagrangians}

\begin{definition}
An \underline{enhanced Lagrangian} is a homomorphism section (with respect
to the canonical projection $H\left( V\right) \rightarrow V$) $\tau
:L\rightarrow H\left( V\right) $, where $L\in Lag\left( V\right) $ is a
Lagrangian subspace in $V$.
\end{definition}

We denote by $ELag\left( V\right) $ the set of enhanced Lagrangians. A
concrete way to specify an enhanced Lagrangian is as follows: Let $L\in
Lag\left( V\right) $, denote by $\Sigma _{L}$ the set consisting of
"quadratic functions" $\alpha :L\rightarrow R$ satisfying 
\begin{equation}
\alpha \left( l_{1}+l_{2}\right) -\alpha \left( l_{1}\right) -\alpha \left(
l_{2}\right) =\beta \left( l_{1},l_{2}\right) \text{.}  \label{enhanced_eq}
\end{equation}

A pair $\left( L,\alpha \right) $ yields a homomorphism section $\tau
:L\rightarrow H\left( V\right) $ given by $\tau \left( l\right) =\left(
l,\alpha \left( l\right) \right) $. Indeed we verify 
\begin{eqnarray*}
\tau \left( l_{1}+l_{2}\right) &=&\left( l_{1}+l_{2},\alpha \left(
l_{1}+l_{2}\right) \right) \\
&=&\left( l_{1}+l_{2},\alpha \left( l_{1}\right) +\alpha \left( l_{2}\right)
+\beta \left( l_{1},l_{2}\right) \right) \\
&=&\tau \left( l_{1}\right) \cdot \tau \left( l_{2}\right) \text{.}
\end{eqnarray*}%
where in the second equality we used the characteristic property of $\alpha $
(Equation (\ref{enhanced_eq})).

There is an evident action of the group $ASp\left( V\right) $ on the set $%
ELag\left( V\right) $ sending an element $\tau :L\rightarrow H\left(
V\right) $ to $Ad_{g}\tau :gL\rightarrow H\left( V\right) $, given by $%
Ad_{g}\tau \left( l\right) =g\tau \left( g^{-1}l\right) .$

We note that the set $\Sigma _{L}$ is a principal homogenous set over the
dual group $L^{\vee }=Hom\left( L,R\right) $. We proceed to explain a
construction of specific elements in $\Sigma _{L}$: Let $\widetilde{L}\in
Lag(\widetilde{V})$ be a free Lagrangian sub-module in $\widetilde{V}$ \
such that $\widetilde{L}/\mathfrak{m}_{K}\widetilde{L}=L$. The cocycle $%
\widetilde{\beta }$ \ yields quadratic form on the module $\widetilde{L}$
which we denote also by $\widetilde{\beta }:\widetilde{L}\rightarrow R$.

\begin{lemma}
\label{tech2_lemma}The quadratic form $\widetilde{\beta }:\widetilde{L}%
\rightarrow R$ factors through $L$ and yields a "quadratic function" $\alpha
_{\widetilde{L}}\in \Sigma _{L}$.
\end{lemma}

For a proof, see Appendix \ref{proofs_sec}.

\textbf{Convention: }To simplify notations, we will often denote an enhanced
Lagrangian $\tau :L\rightarrow H\left( V\right) $ simply by $L$.

\subsubsection{Models associated with enhanced Lagrangians}

We associate with each enhanced Lagrangian $L$, a model $\left( \pi
_{L},H\left( V\right) ,\mathcal{H}_{L}\right) $ of the Heisenberg
representation, as follows: The vector space $\mathcal{H}_{L}$ consists of
functions $f:H\left( V\right) \rightarrow 
\mathbb{C}
$ satisfying $f\left( z\cdot \tau \left( l\right) \cdot h\right) =\psi
\left( z\right) f\left( h\right) $, for every $z\in Z$ and $l\in L$ and the
action $\pi _{L}:H\left( V\right) \rightarrow GL\left( \mathcal{H}%
_{L}\right) $ is given by right translations, namely $\pi _{L}\left(
h\right) \left[ f\right] \left( h^{\prime }\right) =f\left( h^{\prime }\cdot
h\right) $, for every $h,h^{\prime }\in H\left( V\right) $.

The collection of models $\left \{ \mathcal{H}_{L}\right \} $ forms a vector
bundle $\mathfrak{H\rightarrow }ELag\left( V\right) $, with fibers $%
\mathfrak{H}_{L}=\mathcal{H}_{L}$. The vector bundle $\mathfrak{H}$ is
equipped with additional structure of an action $\pi _{L}$ of $H\left(
V\right) $ on each fiber. This suggests the following terminology:

\begin{definition}
Let $n\in 
\mathbb{N}
$. An $H\left( V\right) ^{n}$-vector bundle on $ELag\left( V\right) $ is a
vector bundle $\mathfrak{E}\rightarrow ELag\left( V\right) $, equipped with
a fiberwise action $\pi _{L}:H\left( V\right) ^{n}\rightarrow GL(\mathfrak{E}%
_{L})$, for every $L\in ELag\left( V\right) $.
\end{definition}

In addition, our $\mathfrak{H}$ is equipped with a natural $ASp\left(
V\right) $-equivariant structure, defined as follows: For every $g\in
ASp\left( V\right) $, let $g^{\ast }\mathfrak{H}$ be the $H\left( V\right) $%
-vector bundle with fibers $g^{\ast }\mathfrak{H}_{L}=\mathcal{H}_{gL}$ and
the $g$-twisted Heisenberg action $\pi _{L}^{g}:H\left( V\right) \rightarrow
GL\left( \mathcal{H}_{gL}\right) $, given by $\pi _{L}^{g}\left( h\right)
=\pi _{gL}\left( g\left( h\right) \right) $. The equivariant structure is
the isomorphisms of $H\left( V\right) $-vector bundles 
\begin{equation}
\theta _{g}:g^{\ast }\mathfrak{H\rightarrow H},  \label{equivariant_eq}
\end{equation}%
which on the level of fibers, sends $f\in \mathcal{H}_{gL}$ to $f\circ g\in 
\mathcal{H}_{L}$.

\subsection{The strong Stone-von Neumann property}

We proceed to formulate a stronger form of the Stone-von Neumann property of
the Heisenberg representation. First, we introduce the following terminology:

\begin{definition}
Let $\mathfrak{E}\rightarrow ELag\left( V\right) $ be an $H\left( V\right)
^{n}$-vector bundle. A \underline{trivialization} of $\mathfrak{E}$ is a
system of intertwining isomorphisms $\{E_{M,L}\in \mathrm{Hom}_{H\left(
V\right) ^{n}}(\mathfrak{E}_{L},\mathfrak{E}_{M}):\left( M,L\right) \in
ELag\left( V\right) ^{2}\}$ satisfying the following multiplicativity
condition%
\begin{equation*}
E_{N,M}\circ E_{M,L}=E_{N,L},
\end{equation*}%
for every $N,M,L\in ELag\left( V\right) $.
\end{definition}

\begin{remark}
Intuitively, a trivialization of a $H\left( V\right) ^{n}$-vector bundle $%
\mathfrak{E}\rightarrow ELag\left( V\right) $ might be thought of as a flat
connection, compatible with the Heisenberg action and admitting a trivial
monodromy.
\end{remark}

\begin{theorem}[The strong S-vN property]
\label{SS-vN_thm}The $H\left( V\right) ^{4}$-vector bundle $\mathfrak{H}%
^{\otimes 4}$ admits a natural trivialization $\{T_{M,L}\}$.
\end{theorem}

For a proof, see Section \ref{S-vN_sec}.

\subsection{The Weil gerbe\label{Weilgerbe_sub}}

For us a gerbe, is a groupoid category which consists of a single
isomorphism class of objects. The \underline{band} of a gerbe $\mathcal{G}$
is the group \textrm{Aut}$\left( \mathfrak{E}\right) $, for any object $%
\mathfrak{E\in }\mathcal{G}$. We proceed to describe a gerbe $\mathcal{W}$
which is canonically associated with the vector bundle $\mathfrak{H}$.

An object in $\mathcal{W}$ is a triple $\left( \mathfrak{E}%
,\{E_{M,L}\},\varphi \right) $, where

\begin{itemize}
\item $\mathfrak{E}$ is an $H\left( V\right) $-vector bundle.

\item $\{E_{M,L}\}$ is a trivialization of $\mathfrak{E}$.

\item $\varphi :\mathfrak{E}\overset{\simeq }{\rightarrow }\mathfrak{H}$ is
an isomorphism of $H\left( V\right) $-vector bundles which satisfies that 
\begin{equation*}
\varphi ^{\otimes 4}:\mathfrak{E}^{\otimes 4}\overset{\simeq }{\rightarrow }%
\mathfrak{H}^{\otimes 4}
\end{equation*}%
is an isomorphism of trivialized $H\left( V\right) ^{4}$-vector bundles.
\end{itemize}

A morphism $f\in \mathrm{Mor}_{\mathcal{W}}\left( \mathfrak{E}_{1},\mathfrak{%
E}_{2}\right) $ is a morphism of trivialized $H\left( V\right) $-vector
bundles satisfying $\varphi _{2}^{\otimes 4}\circ f^{\otimes 4}=\varphi
_{1}^{\otimes 4}$.

\begin{proposition}
\label{groupoid_prop}The category $\mathcal{W}$ is a gerbe with band $\mu
_{4}$.
\end{proposition}

For a proof, see appendix \ref{proofs_sec}.

\subsubsection{Action of the affine symplectic group on the Weil gerbe}

There is a natural action of the affine symplectic group $ASp\left( V\right) 
$ on the Weil gerbe $\mathcal{W}$, defined as follows: Every element $g\in
ASp\left( V\right) $ acts by a pull-back functor $g^{\ast }:\mathcal{%
W\rightarrow W}$, sending an object $\left( \mathfrak{E},\{E_{M,L}\},\varphi
\right) $ to the object $\left( g^{\ast }\mathfrak{E},\{g^{\ast
}E_{M,L})\},\varphi ^{g}\right) $, where

\begin{itemize}
\item $g^{\ast }\mathfrak{E}$ is the pull-back vector bundle, equipped with
the $g$-twisted $H\left( V\right) $ action $\pi _{L}^{g}:H\left( V\right)
\rightarrow g^{\ast }\mathfrak{E}_{L}=\mathfrak{E}_{gL}$ given by 
\begin{equation*}
\pi _{L}^{g}\left( h\right) =\pi _{L}\left( gh\right) .
\end{equation*}

\item $\{g^{\ast }E_{M,L}\}$ is the pull-back trivialization given by $%
g^{\ast }E_{M,L}=E_{gM,gL}$.

\item $\varphi ^{g}$ is the isomorphism given by $\varphi ^{g}=\theta
_{g}\circ g^{\ast }\varphi $, where $\theta _{g}:g^{\mathfrak{\ast }}%
\mathfrak{H}\rightarrow \mathfrak{H}$ is the equivariant structure on $%
\mathfrak{H}$ (see Equation (\ref{equivariant_eq})).
\end{itemize}

There is a group $AMp\left( V\right) $ which is naturally associated with
the action of $ASp\left( V\right) $ on $\mathcal{W}$. The group $AMp\left(
V\right) $ consists of pairs $\left( g,\iota \right) $, where $g\in
ASp\left( V\right) $ and $\iota $ is an isomorphism of functors $\iota
:g^{\ast }\overset{\simeq }{\rightarrow }Id$, with $Id$ denoting the
identity functor. The multiplication rule is defined as follows: For $\left(
g,\iota _{g}\right) $, $\left( h,\iota _{h}\right) \in AMp\left( V\right) $,
where $\iota _{g}:g^{\ast }\overset{\simeq }{\rightarrow }Id$ and $\iota
_{h}:h^{\ast }\overset{\simeq }{\rightarrow }Id$, their multiplication $%
\left( g,\iota _{g}\right) \cdot \left( h,\iota _{h}\right) $ is the pair $%
\left( g\cdot h,\iota \right) $, where $\iota :\left( gh\right) ^{\ast }%
\overset{\simeq }{\rightarrow }Id$ is the composition $\iota =\iota
_{h}\circ h^{\ast }\left( \iota _{g}\right) $.

\begin{proposition}
\label{action_prop}The group $AMp\left( V\right) $ is a central extension of
the group $ASp\left( V\right) $ by $\mu _{4}$, in particular it fits into an
exact sequence of groups%
\begin{equation*}
1\rightarrow \mu _{4}\rightarrow AMp\left( V\right) \rightarrow ASp\left(
V\right) \rightarrow 1\text{.}
\end{equation*}
\end{proposition}

For a proof, see appendix \ref{proofs_sec}.

\subsection{The canonical vector space}

Let us denote by $\mathsf{Vect}$ the category of complex vector spaces.
There is a natural fiber functor $\Gamma :\mathcal{W}\rightarrow \mathsf{Vect%
}$, sending $\mathfrak{E}\in $ $\mathcal{W}$ to the vector space of
"horizontal sections" 
\begin{equation*}
\Gamma \left( \mathfrak{E}\right) =\Gamma _{hor}\left( ELag\left( V\right) ,%
\mathfrak{E}\right) ,
\end{equation*}%
which consists of compatible systems of vectors $\left( f_{L}\in \mathfrak{E}%
_{L}:L\in ELag\left( V\right) \right) $ such that $E_{M,L}\left(
f_{L}\right) =f_{M}$, for every $\left( M,L\right) \in ELag\left( V\right)
^{2}$.

There exists a natural homomorphism $\rho :AMp\left( V\right) \rightarrow
Aut\left( \Gamma \right) $ defined as follows: Given an element $\left(
g,\iota \right) \in ASp\left( V\right) $, the automorphism $\rho \left(
g,\iota \right) :\Gamma \rightarrow \Gamma $ is the composition 
\begin{equation*}
\Gamma \overset{\simeq }{\rightarrow }\Gamma \circ g^{\ast }\overset{\Gamma
\left( \iota \right) }{\longrightarrow }\Gamma \circ Id=\Gamma ,
\end{equation*}%
where the first morphism is the tautological isomorphism.

We refer to the homomorphism $\rho $ as the canonical model of the Weil
representation; in more scientific terms, $\rho $ is a representation of the
central extension $AMp\left( V\right) $ on a vector space twisted by the
gerbe $\mathcal{W}$.

We obtain a more traditional realization of the Weil representation, after
choosing a specific object $\mathfrak{E}\in \mathcal{W}$. This yields a
homomorphism 
\begin{equation*}
\rho _{\mathfrak{E}}:AMp\left( V\right) \rightarrow GL\left( \Gamma \left( 
\mathfrak{E}\right) \right) \text{.}
\end{equation*}

As a consequence we proved Theorem \ref{Weilrep_thm}.

\section{The strong Stone-von Neumann property \label{S-vN_sec}}

In this section we describe the construction of the canonical trivialization
of the vector bundle $\mathfrak{H}^{\otimes 4}$, which is asserted in
Theorem \ref{SS-vN_thm}.

\subsection{Canonical intertwining morphisms}

The $H\left( V\right) $-vector bundle $\mathfrak{H}$ admits a partial
connective structure which we are going to describe below.

Let us denote by $U_{2}\subset ELag\left( V\right) ^{2}$ the subset
consisting of pairs of enhanced Lagrangians $\left( M,L\right) $ which are
in general position, that is $M+L=V$. For every $\left( M,L\right) \in U_{2}$%
, there exists a canonical intertwining morphism $F_{M,L}\in \mathrm{Hom}%
_{H\left( V\right) }\left( \mathcal{H}_{L},\mathcal{H}_{M}\right) $, given
by the following averaging formula%
\begin{equation*}
F_{M,L}\left[ f\right] \left( h\right) =\tsum \limits_{m\in M}f\left( \tau
\left( m\right) \cdot h\right) ,
\end{equation*}%
for every $f\in \mathcal{H}_{L}$.

Let us denote by $U_{3}\subset ELag\left( V\right) ^{3}$, the subset
consisting of triples of enhanced Lagrangians $\left( N,M,L\right) $ which
are in general position pairwisely. For every $\left( N,M,L\right) \in U_{3}$%
, we can form two intertwining morphisms in $\mathrm{Hom}_{H\left( V\right)
}\left( \mathcal{H}_{L},\mathcal{H}_{N}\right) $. The first is $F_{N,L}$ and
the second is the composition $F_{N,M}\circ F_{M,L}$. Since $\mathcal{H}_{L}$
and $\mathcal{H}_{M}$ are both isomorphic to the Heisenberg representation,
which is irreducible, we have that $F_{N,L}$ and $F_{N,M}\circ F_{M,L}$ are
proportional. Let us denote by $C\left( N,M,L\right) $ the proportion
constant, that is%
\begin{equation*}
F_{N,M}\circ F_{M,L}=C\left( N,M,L\right) \cdot F_{N,L}\text{.}
\end{equation*}

The function $C:U_{3}\rightarrow 
\mathbb{C}
$, which sends a triple $\left( N,M,L\right) $ to $C\left( N,M,L\right) $ is
a cocycle (with respect to an appropriately defined differential), moreover,
it can be described explicitly. For this we need to introduce some
additional terminology.

Let $r^{L}:M\rightarrow N$ denote the linear map characterized by the
condition 
\begin{equation*}
r^{L}\left( m\right) -m\in L,
\end{equation*}%
for every $m\in M$. Equivalently, $r^{L}$ is characterized by the condition
that 
\begin{equation*}
\omega \left( r^{L}\left( m\right) ,l\right) =\omega \left( m,l\right) ,
\end{equation*}
for every $l\in L$. Let us write $N=\left( N,\alpha _{N}\right) ,M=\left(
M,\alpha _{M}\right) $ and $L=\left( L,\alpha _{L}\right) $. In addition,
let us denote by $Q_{\left( N,M,L\right) }:M\rightarrow R$ the "quadratic
function" given by 
\begin{equation*}
Q_{\left( N,M,L\right) }\left( m\right) =\alpha _{M}\left( m\right) +\alpha
_{N}\left( -r^{L}\left( m\right) \right) -\alpha _{L}\left( m-r^{L}\left(
m\right) \right) -\beta \left( m,r^{L}\left( m\right) \right) .
\end{equation*}

\begin{proposition}
\label{explicit_prop}We have%
\begin{equation*}
C\left( N,M,L\right) =\tsum \limits_{m\in M}\psi \left( Q_{\left(
N,M,L\right) }\left( m\right) \right) .
\end{equation*}
\end{proposition}

For a proof, see Appendix \ref{proofs_sec}.

\begin{theorem}
\label{power_thm}For every $\left( N,M,L\right) \in U_{3}$%
\begin{equation*}
C\left( N,M,L\right) ^{4}=\left( -1\right) ^{d\cdot n}\cdot \left \vert
M\right \vert ^{2}\text{.}
\end{equation*}
\end{theorem}

Theorem \ref{power_thm} will follow from Theorem \ref{Gauss_thm}, which
appear below.

Granting the validity of Theorem \ref{power_thm}, we can exhibit the
canonical trivialization of the vector bundle $\mathfrak{H}^{\otimes 4}$.

First, we note that $\mathfrak{H}^{\otimes 4}$ admits the following partial
trivialization: For every $\left( M,L\right) \in U_{2}$ consider the
intertwining morphism 
\begin{equation*}
T_{M,L}=A_{M,L}\cdot F_{M,L}^{\otimes 4},
\end{equation*}%
where $A_{M,L}$ is the normalization coefficient $A_{M,L}=\frac{\left(
-1\right) ^{d\cdot n}}{\left \vert M\right \vert ^{2}}$.

Evidently, the normalization coefficients $\left \{ A_{M,L}:\left(
M,L\right) \in U_{2}\right \} $ satisfy%
\begin{equation*}
A_{N,M}\cdot A_{M,L}=\frac{\left( -1\right) ^{d\cdot n}}{\left \vert M\right
\vert ^{2}}A_{N,L},
\end{equation*}%
for every $\left( N,M,L\right) \in U_{3}$. Hence, we conclude that 
\begin{equation}
T_{N,M}\circ T_{M,L}=T_{N,L},  \label{mult_eq}
\end{equation}%
for every $\left( N,M,L\right) \in U_{3}$.

\begin{theorem}
\label{trivialization_thm}The partial trivialization $\left \{
T_{M,L}:\left( M,L\right) \in U_{2}\right \} $ extends, in a unique manner,
to a trivialization of $\mathfrak{H}^{\otimes 4}$.
\end{theorem}

For a proof, see Appendix \ref{proofs_sec}.

The rest of this section is devoted to the proof of Theorem \ref{power_thm}.
The proof consists of two main steps. In the first step we message the
formula of the cocycle $C\left( N,M,L\right) $ to fit in the setting of
symmetric spaces over the ring $%
\mathbb{Z}
/4%
\mathbb{Z}
$. In the second step we develop the structure theory of symmetric spaces
over $%
\mathbb{Z}
/4%
\mathbb{Z}
$ \ which is then applied to prove the theorem.

\subsection{Simplification of the cocycle}

Choose lifts $\widetilde{N},\widetilde{M},\widetilde{L}\in Lag(\widetilde{V}%
) $, lying over $N,M$ and $L$ respectively.

We can write 
\begin{eqnarray*}
\alpha _{N} &=&\alpha _{\widetilde{N}}+\sigma _{N}, \\
\alpha _{M} &=&\alpha _{\widetilde{M}}+\sigma _{M}, \\
\alpha _{L} &=&\alpha _{\widetilde{L}}+\sigma _{L},
\end{eqnarray*}%
where $\alpha _{\widetilde{N}}\in \Sigma _{N},\alpha _{\widetilde{M}}\in
\Sigma _{M}$ and $\alpha _{\widetilde{L}}\in \Sigma _{L}$ are the enhanced
structures associated with the liftings (see Lemma \ref{tech2_lemma}) and $%
\sigma _{N}\in N^{\vee },\sigma _{M}\in M^{\vee }$ and $\sigma _{L}\in
L^{\vee }$ are characters (taking values in $R$).

Let us denote $Q=Q_{\left( N,M,L\right) }$. Simple verification reveals that 
\begin{equation*}
Q\left( m\right) =\widetilde{\omega }(r^{\widetilde{L}}\left( \widetilde{m}%
\right) ,\widetilde{m})+\sigma \left( m\right) ,
\end{equation*}%
where $\widetilde{m}$ is any element in $\widetilde{M}$ lying over $m\in M$
and $\sigma $ is the element in $M^{\vee }$ given by $\sigma \left( m\right)
=\sigma _{M}\left( m\right) +\sigma _{N}\left( r^{L}\left( m\right) \right)
-\sigma _{L}\left( m-r^{L}\left( m\right) \right) $.

Let us denote by $\widetilde{\omega }_{\widetilde{L}}:\widetilde{M}\times 
\widetilde{M}\rightarrow R$ the bilinear form $\widetilde{\omega }(r^{%
\widetilde{L}}(\cdot ),\cdot )$ and by $\omega _{L}:M\times M\rightarrow R$
the bilinear form \ $\omega (r^{L}(\cdot ),\cdot )$. Evidently, the form $2$ 
$\widetilde{\omega }_{\widetilde{L}}$ reduces to the form $\omega _{L}$.

Since $\omega _{L}$ is non-degenerate, there exists a unique element $%
m_{\sigma }\in M$ such that $\omega _{L}\left( m_{\sigma },\cdot \right)
=\sigma \left( \cdot \right) $. Choosing an element $\widetilde{m}_{\sigma
}\in \widetilde{M}$ lying over $m_{\sigma }$, we can write%
\begin{equation*}
Q\left( m\right) =\widetilde{\omega }_{\widetilde{L}}(\widetilde{m}+%
\widetilde{m}_{\sigma },\widetilde{m}+\widetilde{m}_{\sigma })-\omega
_{L}\left( m_{\sigma },m_{\sigma }\right) \text{.}
\end{equation*}

Hence%
\begin{equation}
C\left( N,M,L\right) =\psi \left( -\omega _{L}\left( m_{\sigma },m_{\sigma
}\right) \right) \tsum \limits_{m\in M}\psi \left( \widetilde{\omega }_{%
\widetilde{L}}(\widetilde{m},\widetilde{m})\right) \text{.}
\label{simple_eq}
\end{equation}

Let $G([\widetilde{M},tr\left( \widetilde{\omega }_{\widetilde{L}}\right) ])$
denote the "Gauss sum" 
\begin{eqnarray*}
G([\widetilde{M},tr\left( \widetilde{\omega }_{\widetilde{L}}\right) ])
&=&\tsum \limits_{m\in M}\psi \left( \widetilde{\omega }_{\widetilde{L}}(%
\widetilde{m},\widetilde{m})\right) \\
&=&\tsum \limits_{m\in M}e^{\frac{2\pi i}{4}tr\left( \widetilde{\omega }_{%
\widetilde{L}}(\widetilde{m},\widetilde{m})\right) }
\end{eqnarray*}

Since $\psi \left( -\omega _{L}\left( m_{\sigma },m_{\sigma }\right) \right)
\in \mu _{2}$, the assertion that $C\left( N,M,L\right) ^{4}=\left(
-1\right) ^{d\cdot n}\cdot \left \vert M\right \vert ^{2}$ will follow from
the following statement

\begin{theorem}
\label{Gauss_thm}We have 
\begin{equation*}
G([\widetilde{M},tr\left( \widetilde{\omega }_{\widetilde{L}}\right)
])^{4}=\left( -1\right) ^{d\cdot n}\cdot \left \vert M\right \vert ^{2}\text{%
.}
\end{equation*}
\end{theorem}

The proof of Theorem \ref{Gauss_thm}, will appear in Subsection \ref%
{Gauss_sub} after we develop some appropriate formalism.

\subsection{Symmetric spaces}

Let $A$ be a ring and let $\mathcal{B}\left( A\right) $ denote the category
of \textit{symmetric spaces} over $A$. An object in $\mathcal{B}\left(
A\right) $ is a pair $\left( V,B\right) $, where $V$ is a free module over $%
A $ and $B:V\times V\rightarrow A$ is a non-degenerate symmetric bilinear
form on $V$. A morphism $f\in $ $\mathrm{Mor}_{\mathcal{B}\left( A\right)
}\left( \left( V_{1},B_{1}\right) ,\left( V_{2},B_{2}\right) \right) $ is a
map of $A $-modules $f:V_{1}\rightarrow V_{2}$, such that $B_{2}\left(
f\left( v_{1}\right) ,f\left( v_{2}\right) \right) =B_{1}\left(
v_{1},v_{2}\right) $, for every $v_{1},v_{2}\in V$.

The category $\mathcal{B}\left( A\right) $ \ has a monoidal structure given
by the operation of direct-sum of symmetric spaces.

Let us denote by $\mathrm{W}\left( A\right) =\left( \mathrm{Iso}\left( 
\mathcal{B}\left( A\right) \right) ,+\right) $ the associated commutative
monoid, whose elements are isomorphism classes of objects in $\mathcal{B}%
\left( A\right) $ and $+$ is the binary operation induced from the monoidal
structure in $\mathcal{B}\left( A\right) $.

\textbf{Notations.} Given an object $\left( V,B\right) \in \mathcal{B}\left(
A\right) $ we will denote by $\left[ V,B\right] $ its isomorphism class in $%
\mathrm{W}\left( A\right) $. More concretely, given a symmetric matrix $M\in 
\mathrm{Mat}_{n\times n}\left( A\right) $, we will denote by $\left[ M\right]
$ the isomorphism class of $\left( A^{n},B_{M}\right) $, where 
\begin{equation*}
B_{M}\left( \overrightarrow{x},\overrightarrow{y}\right) =\overrightarrow{x}%
^{t}\cdot M\cdot \overrightarrow{y}.
\end{equation*}

\subsubsection{The Discriminant}

There is a basic morphism of monoids $d:\mathrm{W}\left( A\right)
\rightarrow $ $A^{\times }/A^{\times 2}$ called the \textit{discriminant}.
Given an element $\left[ V,B\right] \in \mathrm{W}\left( A\right) $, the
discriminant $d\left( \left[ V,B\right] \right) $ can be defined as follows:
Choose an isomorphism $f:V\simeq A^{n}$, $n=rk\left( V\right) $; let $B_{0}$
denote the standard symmetric form on $A^{n}$ given by 
\begin{equation*}
B_{0}\left( \overrightarrow{x},\overrightarrow{y}\right) =x_{1}\cdot
y_{1}+..+x_{n}\cdot y_{n}\text{.}
\end{equation*}

Define $d\left( \left[ V,B\right] \right) =\det \left( B/f^{\ast }\left(
B_{0}\right) \right) $. This procedure yields an element in $A^{\times
}/A^{\times 2}$ that does not depend on the choice of the isomorphism $f$.

\textbf{Notations: }\ Denote $\mathcal{B=B}\left( 
\mathbb{Z}
/4%
\mathbb{Z}
\right) $ and $\mathrm{W}=\mathrm{W}\left( 
\mathbb{Z}
/4%
\mathbb{Z}
\right) $.

We refer to $\mathrm{W}$ as the \textit{Witt-Grothendieck monoid}.

\subsubsection{The structure of the Witt-Grothendieck monoid}

The fine structure of the Witt-Grothendieck monoid is specified in the
following two propositions. The first proposition asserts that every element
in $\mathrm{W}$ can be written in a standard form as a combination of four
types of generators. The second proposition specifies some basic relations
satisfied by these generators.

\begin{proposition}
\label{structure1_prop}Let $\left[ V,B\right] \in \mathrm{W}$ then 
\begin{equation*}
\left[ V,B\right] =n_{1}\cdot \left[ 1\right] +n_{2}\cdot \lbrack
-1]+n_{3}\cdot 
\begin{bmatrix}
0 & 1 \\ 
1 & 0%
\end{bmatrix}%
+n_{4}\cdot 
\begin{bmatrix}
2 & 1 \\ 
1 & 2%
\end{bmatrix}%
,
\end{equation*}%
where $n_{i}\in 
\mathbb{N}
$, $i=1,2,3,4$.
\end{proposition}

For a proof, see Appendix \ref{proofs_sec}.

\begin{proposition}
\label{structure2_prop}The following relations hold in $\mathrm{W}$
\end{proposition}

\begin{eqnarray}
\begin{bmatrix}
2 & 1 \\ 
1 & 2%
\end{bmatrix}%
+%
\begin{bmatrix}
2 & 1 \\ 
1 & 2%
\end{bmatrix}
&=&%
\begin{bmatrix}
0 & 1 \\ 
1 & 0%
\end{bmatrix}%
+%
\begin{bmatrix}
0 & 1 \\ 
1 & 0%
\end{bmatrix}
\label{rel1_eq} \\
3\cdot \left[ 1\right] &=&\left[ -1\right] +%
\begin{bmatrix}
2 & 1 \\ 
1 & 2%
\end{bmatrix}%
,  \label{rel2_eq} \\
3\cdot \left[ -1\right] &=&\left[ 1\right] +%
\begin{bmatrix}
2 & 1 \\ 
1 & 2%
\end{bmatrix}%
.  \label{rel3_eq}
\end{eqnarray}

For a proof, see Appendix \ref{proofs_sec}.

\subsection{The Gauss character and the Witt group}

We describe a morphism of monoids $G=G_{\psi }:\mathrm{W}\rightarrow 
\mathbb{C}
^{\times }$, which we call the \textit{Gauss character}. It is defined as
follows:

For $\left[ V,B\right] \in \mathrm{W}$, let 
\begin{equation*}
G\left( \left[ V,B\right] \right) =\tsum \limits_{v\in V/2V}\psi \left(
B\left( v,v\right) \right) \text{,}
\end{equation*}%
where we note that the value $G\left( \left[ V,B\right] \right) $ does not
depend on the representative $\left( V,B\right) \in \mathcal{B}$ and summing
over the quotient $V/2V$ makes sense since the quadratic function $B\left(
v,v\right) $ factors through $V/2V$. The morphism $G:\mathrm{W}\rightarrow 
\mathbb{C}
^{\times }$ is a morphism of monoids 
\begin{equation*}
G\left( \left[ V_{1},B_{1}\right] +\left[ V_{2},B_{2}\right] \right)
=G\left( \left[ V_{1},B_{1}\right] \right) \cdot G\left( \left[ V_{2},B_{2}%
\right] \right) ,
\end{equation*}%
for every $\left[ V_{1},B_{1}\right] ,\left[ V_{2},B_{2}\right] \in \mathrm{W%
}$.

\begin{theorem}[Purity theorem]
\label{purity_thm}For every $\left[ V,B\right] \in \mathrm{W}$%
\begin{equation*}
\left \vert G\left( \left[ V,B\right] \right) \right \vert =2^{rk\left(
V\right) /2}=\left \vert V/2V\right \vert ^{1/2}\text{.}
\end{equation*}
\end{theorem}

For a proof, see Appendix \ref{proofs_sec}.

Let $I\subset \mathrm{W}$ be the submonoid consisting of elements $\left[ V,B%
\right] \in \mathrm{W}$ such that $G\left( \left[ V,B\right] \right) \in 2^{%
\mathbb{Z}
}$, we denote $\mathrm{GW}=\mathrm{W}/I$.

\begin{theorem}
\label{GWitt_thm}The submonoid $I$ is generated by the elements 
\begin{equation*}
\left[ 1\right] +\left[ -1\right] ,%
\begin{bmatrix}
0 & 1 \\ 
1 & 0%
\end{bmatrix}%
.
\end{equation*}

Moreover, the quotient monoid $\mathrm{GW}=\mathrm{W}/I$ is a group,
isomorphic to $%
\mathbb{Z}
/8%
\mathbb{Z}
$.
\end{theorem}

For a proof, see Appendix \ref{proofs_sec}.

We refer to the group $\mathrm{GW}$ as the \textit{Witt group}.

\subsection{Proof of Theorem \protect \ref{Gauss_thm}\label{Gauss_sub}}

We are now ready to prove Theorem \ref{Gauss_thm}$.$

We know that $\mathrm{GW}\simeq 
\mathbb{Z}
/8%
\mathbb{Z}
$, let us denote by $\mathbf{4}\in \mathrm{GW}$ the unique non-trivial
element of order $2$.

\begin{proposition}
\label{structure3_prop}For every $\left[ V,B\right] \in \mathrm{GW}$%
\begin{equation*}
4\cdot \left[ V,B\right] =rk\left( V\right) \cdot \mathbf{4}\text{.}
\end{equation*}
\end{proposition}

For a proof, see Appendix \ref{proofs_sec}.

The following statement is a direct consequence of the purity theorem
(Theorem \ref{purity_thm}) and Proposition \ref{structure3_prop}

\begin{corollary}
\label{Gauss_cor}For every $\left[ V,B\right] \in \mathrm{W}$%
\begin{equation*}
G\left( 4\cdot \left[ V,B\right] \right) =\left( -1\right) ^{rk\left(
V\right) }\cdot \left \vert V/2V\right \vert ^{2}\text{.}
\end{equation*}
\end{corollary}

Since $\widetilde{\omega }_{\widetilde{L}}:\widetilde{M}\times \widetilde{M}%
\rightarrow R$ is a non-degenerate symmetric form over the ring $R$, it
implies that $tr(\widetilde{\omega }_{\widetilde{L}}):\widetilde{M}\times 
\widetilde{M}\rightarrow 
\mathbb{Z}
/4%
\mathbb{Z}
$ is a non-degenerate symmetric form over $%
\mathbb{Z}
/4%
\mathbb{Z}
$. Considering $(\widetilde{M},tr(\widetilde{\omega }_{\widetilde{L}}))$ as
an object in $\mathcal{B}$ and denoting by $[\widetilde{M},tr(\widetilde{%
\omega }_{\widetilde{L}})]$ its class in $\mathrm{W}$, we can write 
\begin{eqnarray*}
G([\widetilde{M},tr(\widetilde{\omega }_{\widetilde{L}})])^{4} &=&G(4\cdot
\lbrack \widetilde{M},tr(\widetilde{\omega }_{\widetilde{L}})]) \\
&=&\left( -1\right) ^{rk(\widetilde{M})}\cdot \left \vert M\right \vert
^{2}=\left( -1\right) ^{d\cdot n}\cdot \left \vert M\right \vert ^{2},
\end{eqnarray*}%
where, in the first equality we used the fact that $G$ is a morphism of
monoids, in the second equality we used Corollary \ref{Gauss_cor} and in the
third equality we used the fact that the rank of $\widetilde{M}$ considered
as a module over $%
\mathbb{Z}
/4%
\mathbb{Z}
$ is $d\cdot n$ where we remind that $d=\left[ K:%
\mathbb{Q}
_{2}\right] $ and $\dim V=2n$.

This concludes the proof of the theorem.

\section{Splitting of the canonical vector space\label{canonical-split_sec}}

\subsection{The pull-back vector bundle}

Considering the symplectic module $(\widetilde{V},\widetilde{\omega })$, we
introduce the following terminology:

\begin{definition}
An \underline{oriented Lagrangian} in $\widetilde{V}$ is a pair $(\widetilde{%
L},o)$, where $\widetilde{L}\in Lag(\widetilde{V})$ is a free Lagrangian
sub-module of $\widetilde{V}$ and $o\in \wedge ^{n}\widetilde{L}$ is an
element such that $R\cdot o=L$.
\end{definition}

Let us denote by $OLag(\widetilde{V})$ the set of oriented Lagrangians in $%
\widetilde{V}$.

\textbf{Convention: }To simplify notations, we will often denote an oriented
Lagrangian $(\widetilde{L},o)$ simply by $\widetilde{L}$.

There is a forgetful map $\pi :OLag(\widetilde{V})\rightarrow ELag\left(
V\right) $, sending an oriented Lagrangian $(\widetilde{L},o)$ to the
enhanced Lagrangian $\ (L,\alpha _{\widetilde{L}})$, where $L=\widetilde{L}/2%
\widetilde{L}$ \ and $\alpha _{\widetilde{L}}\in \Sigma _{L}$ is the
enhanced structure associated with the lift $\widetilde{L}\rightarrow L$
(see Lemma \ref{tech2_lemma}).

Let us denote by $\widetilde{\mathfrak{H}}$ the $H\left( V\right) $-vector
bundle given by the pull-back $\widetilde{\mathfrak{H}}=\pi ^{\ast }%
\mathfrak{H}$. The $H\left( V\right) $-vector bundle $\widetilde{\mathfrak{H}%
}$ is equipped with a natural $Sp(\widetilde{V})$-equivariant structure
which is induced from the $ASp\left( V\right) $-equivariant structure on $%
\mathfrak{H}$. In addition, the $H\left( V\right) ^{4}$-vector bundle $%
\widetilde{\mathfrak{H}}^{\otimes 4}$ admits a trivialization $\{T_{%
\widetilde{M},\widetilde{L}}\}$ which is induced from the trivialization of $%
\mathfrak{H}^{\otimes 4}$ (see Theorem \ref{SS-vN_thm}).

\subsection{Square root of the canonical trivialization}

The following statements concerns the existence of a natural square root of
the trivialized $H\left( V\right) ^{4}$-vector bundle $\widetilde{\mathfrak{H%
}}^{\otimes 4}$.

\begin{theorem}[The strong S-vN property - split form]
\label{S-vN-split_thm}The $H\left( V\right) ^{2}$-vector bundle $\widetilde{%
\mathfrak{H}}^{\otimes 2}$ admits a natural trivialization $\{S_{\widetilde{M%
},\widetilde{L}}\}$ which satisfies 
\begin{equation*}
S_{\widetilde{M},\widetilde{L}}^{\otimes 2}=T_{\widetilde{M},\widetilde{L}},
\end{equation*}%
for every $(\widetilde{M},\widetilde{L})\in OLag(\widetilde{V})^{2}$.
\end{theorem}

For a proof, see Section \ref{S-vN-split_sec}.

\subsection{Splitting of the Weil gerbe}

Our goal is to describe a natural splitting of the Weil gerbe $\mathcal{W}$.
In more precise terms, restricting the action of $ASp\left( V\right) $ on
the Weil gerbe $\mathcal{W}$ to an action of the symplectic group $Sp(%
\widetilde{V})$, we construct a gerbe $\mathcal{W}^{s}$, with band $\mu _{2}$%
, equipped with an action of $Sp(\widetilde{V})$; and a faithful functor $S:%
\mathcal{W}^{s}\rightarrow \mathcal{W}$ which is compatible with the $Sp(%
\widetilde{V})$ actions on both sides.

The definition of the gerbe $\mathcal{W}^{s}$ proceeds as follows: There are
two gerbes which are naturally associated with the vector bundle $\widetilde{%
\mathfrak{H}}$. The first gerbe, which we denote by $\widetilde{\mathcal{W}}$%
, \ has band $\mu _{4}$ and it is associated with the trivialization of $%
\widetilde{\mathfrak{H}}^{\otimes 4}$. The second gerbe, which we denote by $%
\widetilde{\mathcal{W}}^{s}$, has band $\mu _{2}$ and it\ is associated with
the trivialization of $\widetilde{\mathfrak{H}}^{\otimes 2}$ (see Theorem %
\ref{S-vN-split_thm}). The definition of these gerbes is in complete analogy
to the definition of the Weil gerbe $\mathcal{W}$ (see Subsection \ref%
{Weilgerbe_sub}).

In addition, there are two evident functors: The first functor is a fully
faithful functor $\pi ^{\ast }:\mathcal{W}\rightarrow \widetilde{\mathcal{W}}
$, sending an object $\mathfrak{E}\in \mathcal{W}$ to its pull-back $\pi
^{\ast }\mathfrak{E}\in \widetilde{\mathcal{W}}$. The second functor is the
obvious faithful functor $\widetilde{S}:\widetilde{\mathcal{W}}%
^{s}\rightarrow \widetilde{\mathcal{W}}$ acting as identity on objects and
morphisms.

We define the gerbe $\mathcal{W}^{s}$ to be the fiber product category%
\begin{equation*}
\mathcal{W}^{s}=\widetilde{\mathcal{W}}^{s}\times _{\widetilde{\mathcal{W}}}%
\mathcal{W}\text{.}
\end{equation*}

In more concrete terms, an object of $\mathcal{W}^{s}$ is a triple $(%
\widetilde{\mathfrak{E}},\mathfrak{E,\alpha )}$, where $\widetilde{\mathfrak{%
E}}\in \widetilde{\mathcal{W}}^{s}$, $\mathfrak{E\in }$ $\mathcal{W}$ and $%
\alpha \in \mathrm{Mor}_{\widetilde{\mathcal{W}}}(\widetilde{S}(\widetilde{%
\mathfrak{E}}),\pi ^{\ast }(\mathfrak{E})\mathfrak{)}$; and a morphism in $%
\mathrm{Mor}_{\mathcal{W}^{s}}((\widetilde{\mathfrak{E}}_{1},\mathfrak{E}_{1}%
\mathfrak{,\alpha }_{1}\mathfrak{)},(\widetilde{\mathfrak{E}}_{2},\mathfrak{E%
}_{2}\mathfrak{,\alpha }_{2}\mathfrak{))}$ is a pair of morphisms $(%
\widetilde{f},f)$, with $\widetilde{f}\in \mathrm{Mor}_{\widetilde{\mathcal{W%
}}^{s}}(\widetilde{\mathfrak{E}}_{1},\widetilde{\mathfrak{E}}_{2}\mathfrak{)}
$ and $f\in \mathrm{Mor}_{\mathcal{W}}(\mathfrak{E}_{1},\mathfrak{E}_{2}%
\mathfrak{)}$ so that the following compatibility condition is satisfied 
\begin{equation*}
\alpha _{2}\circ \widetilde{S}(\widetilde{f})=\pi ^{\ast }\left( f\right)
\circ \alpha _{1}\text{.}
\end{equation*}

\begin{proposition}
\label{gerbe-split_prop}The category $\mathcal{W}^{s}$ is a gerbe with band $%
\mu _{2}$.
\end{proposition}

For a proof, see Appendix \ref{proofs_sec}.

Finally, we define the splitting functor $S:\mathcal{W}^{s}\rightarrow 
\mathcal{W}$ to be the functor which sends an object $(\widetilde{\mathfrak{E%
}},\mathfrak{E,\alpha )\in }\mathcal{W}^{s}$ to the object $\mathfrak{E\in }%
\mathcal{W}$.

\subsubsection{Action of the symplectic group}

The symplectic group $Sp(\widetilde{V})$ naturally acts on the gerbes $%
\mathcal{W},\widetilde{\mathcal{W}},\widetilde{\mathcal{W}}^{s}$ and $%
\mathcal{W}^{s}$. The definition of these actions is in complete analogy to
the definition of the action of the affine symplectic group $ASp\left(
V\right) $ on the Weil gerbe $\mathcal{W}$ (see Subsection \ref%
{Weilgerbe_sub}). Moreover, the functors $\pi ^{\ast },\widetilde{S}$ and $S$
are compatible with the above actions.

In complete analogy to the definition of the central extension $AMp\left(
V\right) $ (see Subsection \ref{Weilgerbe_sub}), there is a central extension%
\begin{equation*}
1\rightarrow \mu _{2}\rightarrow Mp(\widetilde{V})\rightarrow Sp(\widetilde{V%
})\rightarrow 1\text{,}
\end{equation*}%
naturally associated with the action of $Sp(\widetilde{V})$ on the of gerbe $%
\mathcal{W}^{s}$. An element of the group $Mp(\widetilde{V})$ is a pair $%
\left( \widetilde{g},\widetilde{\iota }\right) $ where $\widetilde{g}\in Sp(%
\widetilde{V})$ and $\widetilde{\iota }$ is an isomorphism of functors $%
\iota :\widetilde{g}^{\ast }\overset{\simeq }{\rightarrow }Id$.

\subsection{The splitting homomorphism}

The splitting functor $S:\mathcal{W}^{s}\rightarrow \mathcal{W}$ induces a
homomorphism%
\begin{equation*}
s:Mp(\widetilde{V})\rightarrow AMp\left( V\right) ,
\end{equation*}
sending an element $\left( \widetilde{g},\widetilde{\iota }\right) \in Mp(%
\widetilde{V})$ to the element $\left( g,\iota \right) \in ASp(V)$ where $%
g=\left( \widetilde{g},\alpha _{\widetilde{g}}\right) $ and $\iota :g^{\ast
}\rightarrow Id$ is the composition%
\begin{equation*}
g^{\ast }\left( S\left( \mathfrak{E}\right) \right) \overset{\simeq }{%
\longrightarrow }S\left( g^{\ast }\left( \mathfrak{E}\right) \right) \overset%
{S\left( \widetilde{\iota }\right) }{\longrightarrow }S\left( \mathfrak{E}%
\right) ,
\end{equation*}%
for every $\mathfrak{E\in }$ $\mathcal{W}^{s}$. This proves, in particular,
Theorem \ref{Weilrep-split_thm}.

\section{The Strong S-vN property - split form\label{S-vN-split_sec}}

In this section we describe the construction of the trivialization of the
vector bundle $\widetilde{\mathfrak{H}}^{\otimes 2}$, which is asserted in
Theorem \ref{S-vN-split_thm}.

\subsection{Canonical intertwining morphisms}

Let us denote by $\widetilde{U}_{2}\subset OLag(\widetilde{V})^{2}$ the
subset consisting of pairs of oriented Lagrangians $(\widetilde{M},%
\widetilde{L})$ which are in general position, that is $\widetilde{M}+%
\widetilde{L}=\widetilde{V}$. For every $(\widetilde{M},\widetilde{L})\in 
\widetilde{U}_{2}$, there is a canonical intertwining morphism $F_{%
\widetilde{M},\widetilde{L}}\in \mathrm{Hom}_{H\left( V\right) }\left( 
\mathcal{H}_{L},\mathcal{H}_{M}\right) $, given by averaging 
\begin{equation*}
F_{\widetilde{M},\widetilde{L}}\left[ f\right] \left( h\right) =\tsum
\limits_{m\in M}f\left( \tau \left( m\right) \cdot h\right) ,
\end{equation*}%
for every $f\in \mathcal{H}_{L}$, where $\tau :M\rightarrow H\left( V\right) 
$ is the enhanced Lagrangian $\pi (\widetilde{M})$.

Let us denote by $\widetilde{U}_{3}\subset OLag(\widetilde{V})^{3}$, the
subset consisting of triples of oriented Lagrangians $(\widetilde{N},%
\widetilde{M},\widetilde{L})$ which are in general position pairwisely. For
every $(\widetilde{N},\widetilde{M},\widetilde{L})\in \widetilde{U}_{3}$ the
intertwining morphisms $F_{\widetilde{N},\widetilde{L}}$ and $F_{\widetilde{N%
},\widetilde{M}}\circ F_{\widetilde{M},\widetilde{L}}$ are proportional%
\begin{equation*}
F_{\widetilde{N},\widetilde{M}}\circ F_{\widetilde{M},\widetilde{L}}=C(%
\widetilde{N},\widetilde{M},\widetilde{L})\cdot F_{\widetilde{N},\widetilde{L%
}}\text{.}
\end{equation*}

The cocycle function $C:\widetilde{U}_{3}\rightarrow 
\mathbb{C}
$ can be described explicitly as follows: Let $r^{\widetilde{L}}:\widetilde{M%
}\rightarrow \widetilde{N}$ denote the linear map characterized by the
condition 
\begin{equation*}
r^{\widetilde{L}}\left( \widetilde{m}\right) -\widetilde{m}\in \widetilde{L},
\end{equation*}%
for every $\widetilde{m}\in \widetilde{M}$. Let $\widetilde{\omega }_{%
\widetilde{L}}:\widetilde{M}\times \widetilde{M}\rightarrow R$ denote the
symmetric form given by $\widetilde{\omega }_{\widetilde{L}}\left( 
\widetilde{m}_{1},\widetilde{m}_{2}\right) =$ $\widetilde{\omega }(r^{%
\widetilde{L}}\left( \widetilde{m}_{1}\right) ,\widetilde{m}_{2}).$

\begin{proposition}
\label{explicit1_prop}For every $(\widetilde{N},\widetilde{M},\widetilde{L}%
)\in \widetilde{U}_{3}$%
\begin{equation*}
C(\widetilde{N},\widetilde{M},\widetilde{L})=G([\widetilde{M},tr(\widetilde{%
\omega }_{\widetilde{L}})]).
\end{equation*}
\end{proposition}

For a proof, see Appendix \ref{proofs_sec}.

\subsection{Normalization coefficients}

For every $(\widetilde{M},\widetilde{L})\in \widetilde{U}_{2}$ define the
normalization coefficient 
\begin{equation*}
A_{\widetilde{M},\widetilde{L}}=G(2[R^{n},tr(B_{\widetilde{M},\widetilde{L}%
})]),
\end{equation*}%
where $B_{\widetilde{M},\widetilde{L}}:R^{n}\times R^{n}\rightarrow R$
denote the symmetric bilinear form 
\begin{equation*}
B_{\widetilde{M},\widetilde{L}}\left( \overrightarrow{x},\overrightarrow{y}%
\right) =x_{1}\cdot y_{1}+..+x_{n-1}\cdot y_{n-1}+\widetilde{\omega }%
_{\wedge }(o_{\widetilde{L}},o_{\widetilde{M}})x_{n}\cdot y_{n},
\end{equation*}%
where $\widetilde{\omega }_{\wedge }:\wedge ^{n}\widetilde{L}\times \wedge
^{n}\widetilde{M}\rightarrow R$ is the pairing induced from the symplectic
form $\widetilde{\omega }$.

\begin{theorem}
\label{cocycle_thm}For every $(\widetilde{N},\widetilde{M},\widetilde{L})\in 
\widetilde{U}_{3}$%
\begin{equation*}
A_{\widetilde{N},\widetilde{M}}\circ A_{\widetilde{M},\widetilde{L}}=G(2[%
\widetilde{M},-tr(\widetilde{\omega }_{\widetilde{L}})])\cdot A_{\widetilde{N%
},\widetilde{L}}.
\end{equation*}
\end{theorem}

The proof of Theorem \ref{cocycle_thm} appears in Subsection \ref%
{cocycle_sub}.

\subsection{Normalized intertwining morphisms}

We are now ready to exhibit the trivialization of the vector bundle $%
\widetilde{\mathfrak{H}}^{\otimes 2}$.

First we define a partial trivialization as follows: Let $S_{\widetilde{M},%
\widetilde{L}}$ be the normalized intertwining morphism 
\begin{equation*}
S_{\widetilde{M},\widetilde{L}}=\frac{A_{\widetilde{M},\widetilde{L}}}{\left
\vert M\right \vert ^{2}}\cdot F_{\widetilde{M},\widetilde{L}}^{\otimes 2},
\end{equation*}%
For every $(\widetilde{M},\widetilde{L})\in \widetilde{U}_{2}$. We have for
every $(\widetilde{N},\widetilde{M},\widetilde{L})\in \widetilde{U}_{3}$

\begin{eqnarray*}
S_{\widetilde{N},\widetilde{M}}\circ S_{\widetilde{M},\widetilde{L}} &=&%
\frac{A_{\widetilde{N},\widetilde{M}}\cdot A_{\widetilde{M},\widetilde{L}}}{%
\left \vert M\right \vert ^{4}}F_{\widetilde{N},\widetilde{L}}^{\otimes
2}\circ F_{\widetilde{M},\widetilde{L}}^{\otimes 2} \\
&=&\frac{G(2[\widetilde{M},-tr(\widetilde{\omega }_{\widetilde{L}})])\cdot
A_{\widetilde{N},\widetilde{L}}}{\left \vert M\right \vert ^{4}}G(2[%
\widetilde{M},tr(\widetilde{\omega }_{\widetilde{L}})])\cdot F_{\widetilde{N}%
,\widetilde{L}}^{\otimes 2} \\
&=&\frac{G(2[\widetilde{M},-tr(\widetilde{\omega }_{\widetilde{L}})]+2[%
\widetilde{M},tr(\widetilde{\omega }_{\widetilde{L}})])}{\left \vert M\right
\vert ^{2}}S_{\widetilde{N},\widetilde{L}}=S_{\widetilde{N},\widetilde{L}},
\end{eqnarray*}%
where in the second equality we used Proposition\  \ref{explicit1_prop} and
Theorem \ref{cocycle_thm} and in the third equality we used the purity
theorem (Theorem \ref{purity_thm}) for the Gauss character and the fact that 
$G\left( \left[ V,B\right] \right) =\overline{G\left( \left[ V,-B\right]
\right) }$, for every $\left( V,B\right) \in \mathcal{B}$.

\begin{theorem}
\label{trivialization-split_thm}The partial trivialization $\{S_{\widetilde{M%
},\widetilde{L}}:(\widetilde{M},\widetilde{L})\in \widetilde{U}_{2}\}$
extends, in a unique manner, to a trivialization of $\widetilde{\mathfrak{H}}%
^{\otimes 2}$.
\end{theorem}

The proof of Theorem \ref{trivialization-split_thm} is similar to the one of
Theorem \ref{trivialization_thm}.

We are left to show that for every $(\widetilde{M},\widetilde{L})\in OLag(%
\widetilde{V})^{2}$ 
\begin{equation*}
S_{\widetilde{M},\widetilde{L}}^{\otimes 2}=T_{\widetilde{M},\widetilde{L}},
\end{equation*}

It is enough to verify this in the case $(\widetilde{M},\widetilde{L})\in 
\widetilde{U}_{2}$, which can be done by direct computation 
\begin{eqnarray*}
S_{\widetilde{M},\widetilde{L}}^{\otimes 2} &=&\frac{G(2[R^{n},tr(B_{%
\widetilde{M},\widetilde{L}})])^{2}}{\left \vert M\right \vert ^{4}}F_{%
\widetilde{M},\widetilde{L}}^{\otimes 2} \\
&=&\frac{G(4[R^{n},tr(B_{\widetilde{M},\widetilde{L}})])}{\left \vert
M\right \vert ^{4}}F_{\widetilde{M},\widetilde{L}}^{\otimes 2} \\
&=&\frac{\left( -1\right) ^{rk\left( R^{n}\right) }\cdot \left \vert M\right
\vert ^{2}}{\left \vert M\right \vert ^{4}}F_{\widetilde{M},\widetilde{L}%
}^{\otimes 2}=\frac{\left( -1\right) ^{d\cdot n}}{\left \vert M\right \vert
^{2}}F_{\widetilde{M},\widetilde{L}}^{\otimes 2}=T_{\widetilde{M},\widetilde{%
L}}\text{,}
\end{eqnarray*}%
where in the second equality we used the fact that $G$ is a morphism of
monoids, in the third equality we used Corollary \ref{Gauss_cor} and in the
forth equality we used the fact that the rank of $R^{n}$ as a module over $%
\mathbb{Z}
/4%
\mathbb{Z}
$ is $d\cdot n$.

\subsection{Proof of Theorem \protect \ref{cocycle_thm}\label{cocycle_sub}}

First we note that by the purity theorem (Theorem \ref{purity_thm}) it is
enough to show that in $\mathrm{GW}$ the following relation holds%
\begin{equation*}
2[\widetilde{M},tr(\widetilde{\omega }_{\widetilde{L}})]+2[tr(B_{\widetilde{N%
},\widetilde{M}})]+2[tr(B_{\widetilde{M},\widetilde{L}})]+2[-tr(B_{%
\widetilde{M},\widetilde{L}})]=\mathbf{0},
\end{equation*}%
where we use the abbreviated notation $[tr(B)]$ for $[R^{n},tr(B)]$.

Let us denote by $X$ the element in $\mathrm{W}$%
\begin{equation*}
X=[\widetilde{M},tr(\widetilde{\omega }_{\widetilde{L}})]+[tr(B_{\widetilde{N%
},\widetilde{M}})]+[tr(B_{\widetilde{M},\widetilde{L}})]+[-tr(B_{\widetilde{M%
},\widetilde{L}})].
\end{equation*}

As an element in $\mathrm{W}$ we have that $rk\left( X\right) =4n$, in
particular, $4|rk\left( X\right) $. Moreover

\begin{proposition}
\label{disc_prop}The discriminant $d\left( X\right) =1.$
\end{proposition}

For a proof, see Appendix \ref{proofs_sec}.

The theorem now follows from

\begin{proposition}
\label{vanishing_prop}Let $X\in \mathrm{W}$ such that $4|rk\left( X\right) $
and $d\left( X\right) =1$ then $2X=\mathbf{0}$ in $\mathrm{GW}$.
\end{proposition}

For a proof, see Appendix \ref{proofs_sec}.

This concludes the proof of the theorem.

\appendix

\section{Proof of statements\label{proofs_sec}}

\subsection{Proof of Lemma \protect \ref{tech1_lemma}}

First we show that the map $\widetilde{\alpha }$ factors to a function $%
\alpha _{\widetilde{g}}:V\rightarrow R$. Let $\widetilde{v},\widetilde{x}\in 
\widetilde{V}$. Write%
\begin{eqnarray*}
\alpha _{\widetilde{g}}\left( \widetilde{v}+2\widetilde{x}\right) &=&\alpha
_{\widetilde{g}}\left( \widetilde{v}\right) +\beta \left( \widetilde{g}%
\left( \widetilde{x}\right) ,\widetilde{g}\left( \widetilde{v}\right)
\right) +\beta \left( \widetilde{g}\left( \widetilde{v}\right) ,\widetilde{g}%
\left( \widetilde{x}\right) \right) -\beta \left( \widetilde{x},\widetilde{v}%
\right) -\beta \left( \widetilde{v},\widetilde{x}\right) \\
&=&\alpha _{\widetilde{g}}\left( \widetilde{v}\right) +\omega \left( 
\widetilde{g}\left( \widetilde{v}\right) ,\widetilde{g}\left( \widetilde{x}%
\right) \right) -\omega \left( \widetilde{v},\widetilde{x}\right) =\alpha _{%
\widetilde{g}}\left( \widetilde{v}\right) .
\end{eqnarray*}

where in the second equality we added $\beta \left( \widetilde{g}\left( 
\widetilde{x}\right) ,\widetilde{g}\left( \widetilde{v}\right) \right)
-\beta \left( \widetilde{g}\left( \widetilde{x}\right) ,\widetilde{g}\left( 
\widetilde{v}\right) \right) $ and also $\beta \left( \widetilde{x},%
\widetilde{v}\right) -\beta \left( \widetilde{x},\widetilde{v}\right) $ and
used the fact that $2\beta =0$.

We are left to show that $\alpha _{\widetilde{g}}\in \Sigma _{g}$. Write 
\begin{eqnarray*}
\alpha _{\widetilde{g}}\left( v_{1}+v_{2}\right) &=&\widetilde{\beta }\left( 
\widetilde{g}\left( \widetilde{v}_{1}+\widetilde{v}_{2}\right) ,\widetilde{g}%
\left( \widetilde{v}_{1}+\widetilde{v}_{2}\right) \right) -\widetilde{\beta }%
\left( \widetilde{v}_{1}+\widetilde{v}_{2},\widetilde{v}_{1}+\widetilde{v}%
_{2}\right) \\
&=&\alpha _{\widetilde{g}}\left( v_{1}\right) +\alpha _{\widetilde{g}}\left(
v_{2}\right) +\widetilde{\beta }\left( \widetilde{g}\left( \widetilde{v}%
_{1}\right) ,\widetilde{g}\left( \widetilde{v}_{2}\right) \right) +%
\widetilde{\beta }\left( \widetilde{g}\left( \widetilde{v}_{2}\right) ,%
\widetilde{g}\left( \widetilde{v}_{1}\right) \right) \\
&&-\widetilde{\beta }\left( \widetilde{v}_{1},\widetilde{v}_{2}\right) -%
\widetilde{\beta }\left( \widetilde{v}_{2},\widetilde{v}_{1}\right) \\
&=&\alpha _{\widetilde{g}}\left( v_{1}\right) +\alpha _{\widetilde{g}}\left(
v_{2}\right) +\beta \left( g\left( v_{1}\right) ,g\left( v_{2}\right)
\right) +\widetilde{\omega }\left( \widetilde{g}\left( \widetilde{v}%
_{2}\right) ,\widetilde{g}\left( \widetilde{v}_{1}\right) \right) \\
&&-\beta \left( v_{1},v_{2}\right) -\widetilde{\omega }\left( \widetilde{v}%
_{2},\widetilde{v}_{1}\right) \\
&=&\alpha _{\widetilde{g}}\left( v_{1}\right) +\alpha _{\widetilde{g}}\left(
v_{2}\right) +\beta \left( g\left( v_{1}\right) ,g\left( v_{2}\right)
\right) -\beta \left( v_{1},v_{2}\right) \text{,}
\end{eqnarray*}%
where, in the third equality we added $\widetilde{\beta }\left( \widetilde{g}%
\left( \widetilde{v}_{1}\right) ,\widetilde{g}\left( \widetilde{v}%
_{2}\right) \right) -\widetilde{\beta }\left( \widetilde{g}\left( \widetilde{%
v}_{1}\right) ,\widetilde{g}\left( \widetilde{v}_{2}\right) \right) $ and
also $\widetilde{\beta }\left( \widetilde{v}_{1},\widetilde{v}_{2}\right) -%
\widetilde{\beta }\left( \widetilde{v}_{1},\widetilde{v}_{2}\right) $ and
used that 
\begin{eqnarray*}
\widetilde{\beta }\left( \widetilde{g}\left( \widetilde{v}_{2}\right) ,%
\widetilde{g}\left( \widetilde{v}_{1}\right) \right) -\widetilde{\beta }%
\left( \widetilde{g}\left( \widetilde{v}_{1}\right) ,\widetilde{g}\left( 
\widetilde{v}_{2}\right) \right) &=&\widetilde{\omega }\left( \widetilde{g}%
\left( \widetilde{v}_{2}\right) ,\widetilde{g}\left( \widetilde{v}%
_{1}\right) \right) , \\
\widetilde{\beta }\left( \widetilde{v}_{2},\widetilde{v}_{1}\right) -%
\widetilde{\beta }\left( \widetilde{v}_{1},\widetilde{v}_{2}\right) &=&%
\widetilde{\omega }\left( \widetilde{v}_{2},\widetilde{v}_{1}\right) \text{.}
\end{eqnarray*}

This concludes the proof of the lemma.

\subsection{Proof of Lemma \protect \ref{tech2_lemma}}

First we show that $\widetilde{\beta }$ factors to a function $\alpha _{%
\widetilde{L}}:L\rightarrow R$. Let $\widetilde{l},\widetilde{x}\in 
\widetilde{L}$. Write%
\begin{eqnarray*}
\widetilde{\beta }(\widetilde{l}+2\widetilde{x},\widetilde{l}+2\widetilde{x}%
) &=&\widetilde{\beta }(\widetilde{l},\widetilde{l})+2\widetilde{\beta }(%
\widetilde{x},\widetilde{l})+2\widetilde{\beta }(\widetilde{l},\widetilde{x})
\\
&=&\widetilde{\beta }(\widetilde{l},\widetilde{l})+4\widetilde{\beta }(%
\widetilde{x},\widetilde{l})\text{,}
\end{eqnarray*}%
where, in the second equality we used the fact that $\widetilde{\beta }:%
\widetilde{L}\times \widetilde{L}\rightarrow R$ is symmetric, since $%
\widetilde{L}$ is a Lagrangian sub-module.

We are left to show that $\alpha _{\widetilde{L}}\in \Sigma _{L}$. Write%
\begin{eqnarray*}
\alpha _{\widetilde{L}}(l_{1}+l_{2})-\alpha _{\widetilde{L}}(l_{1})-\alpha _{%
\widetilde{L}}(l_{2}) &=&\widetilde{\beta }(\widetilde{l}_{1}+\widetilde{l}%
_{2},\widetilde{l}_{1}+\widetilde{l}_{2})-\widetilde{\beta }(\widetilde{l}%
_{1},\widetilde{l}_{1})-\widetilde{\beta }(\widetilde{l}_{2},\widetilde{l}%
_{2}) \\
&=&\widetilde{\beta }(\widetilde{l}_{1},\widetilde{l}_{2})+\widetilde{\beta }%
(\widetilde{l}_{2},\widetilde{l}_{1})=2\widetilde{\beta }(\widetilde{l}_{1},%
\widetilde{l}_{2}) \\
&=&\beta (l_{1},l_{2}),
\end{eqnarray*}%
where in the third equality we, again, used the fact that $\widetilde{\beta }%
:\widetilde{L}\times \widetilde{L}\rightarrow R$ is symmetric. This
concludes the proof of the lemma.

\subsection{Proof of Proposition \protect \ref{groupoid_prop}}

First we show that $\mathcal{W}$ is a groupoid. Let $f:\mathfrak{E}%
_{1}\rightarrow \mathfrak{E}_{2}$, where $\mathfrak{E}_{1},\mathfrak{E}%
_{2}\in \mathcal{W}$. Since $\varphi _{2}^{\otimes 4}\circ f^{\otimes
4}=\varphi _{1}^{\otimes 4}$ and $\varphi _{2}^{\otimes 4}$, $\varphi
_{1}^{\otimes 4}$ are isomorphisms, this implies that $f^{\otimes 4}$ is an
isomorphism, which, in turns, implies that $f$ is an isomorphism.

Second, we show that every two objects in $\mathcal{W}$ are isomorphic. Let $%
\mathfrak{E}_{1},\mathfrak{E}_{2}\in \mathcal{W}$, in order to specify an
isomorphism between $\mathfrak{E}_{1}$ and $\mathfrak{E}_{2}$ as $H\left(
V\right) $-vector bundles with trivializations, it is enough to specify an
isomorphism of a single fiber. Let $L\in ELag\left( V\right) $, and choose $%
f_{L}\in \mathrm{Hom}_{H\left( V\right) }(\mathfrak{E}_{1,L},\mathfrak{E}%
_{2,L})$ to be a non-zero intertwining isomorphism which satisfies $\varphi
_{2,L}^{\otimes 4}\circ f_{L}^{\otimes 4}=\varphi _{1,L}^{\otimes 4}$ (this
can be always done). This implies that the isomorphism $f:\mathfrak{E}%
_{1}\rightarrow \mathfrak{E}_{2}$ which is determined by $f_{L}$ satisfies $%
\varphi _{2}^{\otimes 4}\circ f^{\otimes 4}=\varphi _{1}^{\otimes 4}$ which
means that $f\in \mathrm{Mor}_{\mathcal{W}}(\mathfrak{E}_{1},\mathfrak{E}%
_{2})$.

Finally, we show that, $\mathrm{Mor}\left( \mathfrak{E},\mathfrak{E}\right)
\simeq \mu _{4}$, for every $\mathfrak{E}$ $\in \mathcal{W}$. Let $\mathfrak{%
E}$ $\in \mathcal{W}$. It is easy to see that any two morphisms $%
f_{1},f_{2}\in \mathrm{Mor}\left( \mathfrak{E},\mathfrak{E}\right) $ are
proportional, namely $f_{2}=\lambda \cdot f_{1}$, for some $\lambda \in 
\mathbb{C}
^{\times }$. Now since both satisfies the condition $\varphi ^{\otimes
4}\circ f_{i}^{\otimes 4}=\varphi ^{\otimes 4}$, this implies that $\lambda
^{4}=1$.

This concludes the proof of the proposition.

\subsection{Proof of Proposition \protect \ref{action_prop}}

First we show that $AMp\left( V\right) $ fits into an exact sequence 
\begin{equation*}
1\rightarrow \mu _{4}\rightarrow AMp\left( V\right) \overset{p}{\rightarrow }%
ASp\left( V\right) \rightarrow 1\text{.}
\end{equation*}%
where the the morphism $p:AMp\left( V\right) \rightarrow ASp\left( V\right) $
is the canonical projection. The kernel of $p$ consists of pairs of the form 
$\left( 1,\iota \right) $, where $\iota :Id\overset{\simeq }{\rightarrow }Id$%
. An isomorphism $\iota :Id\overset{\simeq }{\rightarrow }Id$ is determined
by $\iota _{\mathfrak{E}}:\mathfrak{E}\overset{\simeq }{\rightarrow }%
\mathfrak{E}$, for any $\mathfrak{E}\in \mathcal{W}$, but by Proposition \ref%
{groupoid_prop} we know that $\mathrm{Mor}\left( \mathfrak{E},\mathfrak{E}%
\right) \simeq \mu _{4}$.

Second we show that $p:AMp\left( V\right) \rightarrow ASp\left( V\right) $
is a central extension. We have to show that 
\begin{equation*}
\left( g,\iota _{g}\right) ^{-1}\cdot \left( 1,\iota \right) \cdot \left(
g,\iota _{g}\right) =\left( 1,\iota \right) ,
\end{equation*}%
for every $\left( g,\iota _{g}\right) \in AMp\left( V\right) $. Explicit
calculation reveals that 
\begin{equation*}
\left( g,\iota _{g}\right) ^{-1}=\left( g^{-1},\left( g^{-1}\right) ^{\ast
}\left( \iota _{g}^{-1}\right) \right) ,
\end{equation*}
which implies that 
\begin{equation*}
\left( g,\iota _{g}\right) ^{-1}\cdot \left( 1,\iota \right) \cdot \left(
g,\iota _{g}\right) =\left( 1,\iota _{g}\circ \iota \circ \iota
_{g}^{-1}\right) \text{.}
\end{equation*}

Now, verify that $\iota _{g}\circ \iota \circ \iota _{g}^{-1}=\iota $. This
concludes the proof of the proposition.

\subsection{Proof of Proposition \protect \ref{explicit_prop}}

Let $\delta \in \mathcal{H}_{L}$ be the unique function on $H\left( V\right) 
$, supported on $Z\cdot L\subset H\left( V\right) $ and normalized such that 
$\delta \left( 0\right) =1$.

On the one hand, explicit computation reveals that $F_{N,L}\left[ \delta %
\right] \left( 0\right) =1$. This implies that 
\begin{equation}
C\left( N,M,L\right) =F_{N,M}\circ F_{M,L}\left[ \delta \right] \left(
0\right) \text{.}  \label{eq1}
\end{equation}

On the other hand%
\begin{equation*}
F_{N,M}\circ F_{M,L}\left[ \delta \right] \left( 0\right) =\sum
\limits_{m\in M}\sum \limits_{n\in N}\delta \left( \tau _{M}\left( m\right)
\cdot \tau _{N}\left( n\right) \right) \text{,}
\end{equation*}%
where $\tau _{M}:M\rightarrow H\left( V\right) $ and $\tau _{N}:N\rightarrow
H\left( V\right) $ are the associated injective homomorphisms. The
multiplication rule in $H\left( V\right) $ implies that for every $m\in M$
and $n\in N$ we have that 
\begin{equation*}
\tau _{M}\left( m\right) \cdot \tau _{N}\left( n\right) =\left( m+n,\alpha
_{M}\left( m\right) +\alpha _{N}\left( n\right) +\beta \left( m,n\right)
\right) .
\end{equation*}

Since $\delta \left( \tau _{M}\left( m\right) \cdot \tau _{N}\left( n\right)
\right) =0$ unless $m+n\in L$ we obtain that for every $m\in M$, the only
non-zero contribution to the sum $\sum \limits_{n\in N}\delta \left( \tau
_{M}\left( m\right) \cdot \tau _{N}\left( n\right) \right) $ comes from $%
n=-r^{L}\left( m\right) $.

Therefore, we obtain that%
\begin{eqnarray*}
F_{N,M}\circ F_{M,L}\left[ \delta \right] \left( 0\right) &=&\sum
\limits_{m\in M}\delta \left( \tau _{M}\left( m\right) \cdot \tau _{N}\left(
-r^{L}\left( m\right) \right) \right) \\
&=&\sum \limits_{m\in M}\delta \left( m-r^{L}\left( m\right) ,\alpha
_{L}\left( m-r^{L}\left( m\right) \right) +Q_{\left( N,M,L\right) }\left(
m\right) \right) \\
&=&\sum \limits_{m\in M}\psi \left( Q_{\left( N,M,L\right) }\left( m\right)
\right) \text{.}
\end{eqnarray*}

Combining with (\ref{eq1}) we get 
\begin{equation*}
C\left( N,M,L\right) =\sum \limits_{m\in M}\psi \left( Q_{\left(
N,M,L\right) }\left( m\right) \right) \text{.}
\end{equation*}

This concludes the proof of the proposition.

\subsection{Proof of Theorem \protect \ref{trivialization_thm}}

The trivialization of the vector bundle $\mathfrak{H}^{\otimes 4}$ is
constructed as follows: Let $\left( N,L\right) \in ELag\left( V\right) ^{2}$%
, choose a third $M\in ELag\left( V\right) $ such that $\left( N,M\right)
,\left( M,L\right) \in U_{2}$ (such a choice always exists). Define 
\begin{equation*}
T_{N,L}=T_{N,M}\circ T_{M,L}\text{,}
\end{equation*}

Noting that both operators in the left hand side are defined. We are left to
show that the operator $T_{N,L}$ does not depends on the choice of $M$.

Let $M_{i}\in ELag\left( V\right) $, $i=1,2,$ such that $\left(
N,M_{i}\right) ,\left( M_{i},L\right) \in U_{2}$. We want to show that $%
T_{N,M_{1}}\circ T_{M_{1},L}=T_{N,M_{2}}\circ T_{M_{2},L}$.

Choose $M_{3}\in ELag\left( V\right) $ such that $\left( M_{3},M_{i}\right)
\in U_{2}$ and $\left( M_{3},L\right) ,\left( M_{3},N\right) \in U_{2}$. We
have%
\begin{eqnarray*}
T_{N,M_{1}}\circ T_{M_{1},L} &=&T_{N,M_{1}}\circ T_{M_{1},M_{3}}\circ
T_{M_{3},L} \\
&=&T_{N,M_{3}}\circ T_{M_{3},L},
\end{eqnarray*}%
where the first and second equalities are the multiplicativity property for
triples which are in general position pairwisely (Formula \ref{mult_eq}). In
the same fashion, we show that $T_{N,M_{1}}\circ
T_{M_{1},L}=T_{N,M_{3}}\circ T_{M_{3},L}$.

The fact the full system $\left \{ T_{M,L}:\left( M,L\right) \in ELag\left(
V\right) ^{2}\right \} $ is a trivialization can be easily proved along the
same lines as above. This concludes the proof of the proposition.

\subsection{Proof of Proposition \protect \ref{structure1_prop}}

The proof is by induction on the rank of $V$.

If $rk\left( V\right) =1$ then, since $B$ is non-degenerate, either $\left[
V,B\right] =[1]$ or $\left[ V,B\right] =\left[ -1\right] $ ($1$ and $-1$ are
the invertible elements in $%
\mathbb{Z}
/4%
\mathbb{Z}
$) and we are done.

Assume $rk\left( V\right) >1$.

\subsubsection{Case 1}

We are in the situation where there exists $v\in V$ such that $B\left(
v,v\right) =\pm 1$, then the module $V$ decomposes into a direct sum%
\begin{equation*}
V=%
\mathbb{Z}
/4%
\mathbb{Z}
\cdot v\oplus \left( 
\mathbb{Z}
/4%
\mathbb{Z}
\cdot v\right) ^{\bot }\text{,}
\end{equation*}%
therefore $\left[ V,B\right] =\left[ \pm 1\right] +\left[ V_{1},B_{1}\right] 
$, with $rk\left( V_{1}\right) <rk\left( V\right) $ and the statement
follows by induction.

\subsubsection{Case 2}

We are in the situation where $B\left( v,v\right) \in 2\cdot 
\mathbb{Z}
/4%
\mathbb{Z}
$, for every $v\in V$. Always there exists a pair $e,f\in V$ such that $%
B\left( e,f\right) =1$. Let us denote by $U$ the sub-module $%
\mathbb{Z}
/4%
\mathbb{Z}
\cdot e+%
\mathbb{Z}
/4%
\mathbb{Z}
\cdot f$. First we show that $V$ decomposes into a direct sum 
\begin{equation*}
V=U\oplus V_{1}\text{.}
\end{equation*}

We verify the last assertion in three different cases.

\textbf{Case 1. }Assume $B\left( e,e\right) =B\left( f,f\right) =0$. Let $%
P:V\rightarrow V$ be the operator defined by 
\begin{equation*}
P\left( v\right) =v-\left \langle v,e\right \rangle \cdot f-\left \langle
v,f\right \rangle \cdot e.
\end{equation*}

Direct verification reveals that $P$ is an idempotent and $P\left( U\right)
=0$, hence $V=U\oplus PV$.

\textbf{Case 2}. Assume $B\left( e,e\right) =B\left( f,f\right) =2$. Let $%
P:V\rightarrow V$ be the operator defined by 
\begin{equation*}
P\left( v\right) =v-\left \langle v,e\right \rangle \cdot f-\left \langle
v,f\right \rangle \cdot e+2\left \langle v,e\right \rangle \cdot e+2\left
\langle v,f\right \rangle \cdot f\text{.}
\end{equation*}

Direct verification reveals that $P$ is an idempotent and $P\left( U\right)
=0$, hence $V=U\oplus PV$.

\textbf{Case 3}. Without loss of generality, assume $B\left( e,e\right) =2$
and $B\left( f,f\right) =0$. Let $P:V\rightarrow V$ be the operator defined
by 
\begin{equation*}
P\left( v\right) =v-\left \langle v,e\right \rangle \cdot f-\left \langle
v,f\right \rangle \cdot e+2\left \langle v,f\right \rangle \cdot f\text{.}
\end{equation*}

Direct verification reveals that $P$ is an idempotent and $P\left( U\right)
=0$, hence $V=U\oplus PV$.

Now, in case 1 
\begin{equation*}
\left[ U,B_{|U}\right] =%
\begin{bmatrix}
0 & 1 \\ 
1 & 0%
\end{bmatrix}%
,
\end{equation*}%
in case 2 
\begin{equation*}
\left[ U,B_{|U}\right] =%
\begin{bmatrix}
2 & 1 \\ 
1 & 2%
\end{bmatrix}%
,
\end{equation*}%
and in case 3 
\begin{equation*}
\left[ U,B_{|U}\right] =%
\begin{bmatrix}
2 & 1 \\ 
1 & 0%
\end{bmatrix}%
.
\end{equation*}

Finally, considering case 3, we claim that 
\begin{equation*}
\begin{bmatrix}
2 & 1 \\ 
1 & 0%
\end{bmatrix}%
=%
\begin{bmatrix}
0 & 1 \\ 
1 & 0%
\end{bmatrix}%
,
\end{equation*}%
via the isomorphism sending $e\mapsto e+f$ and $f\mapsto f$.

This concludes the proof of the proposition.

\subsection{Proof of Proposition \protect \ref{structure2_prop}}

First we show that 
\begin{equation*}
\begin{bmatrix}
2 & 1 \\ 
1 & 2%
\end{bmatrix}%
+%
\begin{bmatrix}
2 & 1 \\ 
1 & 2%
\end{bmatrix}%
=%
\begin{bmatrix}
0 & 1 \\ 
1 & 0%
\end{bmatrix}%
+%
\begin{bmatrix}
0 & 1 \\ 
1 & 0%
\end{bmatrix}%
\text{.}
\end{equation*}

Denote 
\begin{equation*}
M=%
\begin{bmatrix}
2 & 1 \\ 
1 & 2%
\end{bmatrix}%
,N=%
\begin{bmatrix}
0 & 1 \\ 
1 & 0%
\end{bmatrix}%
\end{equation*}

Let $V_{i}=\left( 
\mathbb{Z}
/4%
\mathbb{Z}
\cdot e_{i}\oplus 
\mathbb{Z}
/4%
\mathbb{Z}
\cdot f_{i},B_{M}\right) $, $i=1,2$ be two copies of the symmetric space
associated to the matrix $M$.

Let $\phi :V_{1}\oplus V_{2}\rightarrow V_{1}\oplus V_{2}$ be the
isomorphism given by 
\begin{eqnarray*}
\phi \left( e_{1}\right) &=&-e_{1}-f_{1}-2e_{2}-f_{2}, \\
\phi \left( f_{1}\right) &=&2e_{1}+f_{1}+2e_{2}-f_{2}, \\
\phi \left( e_{2}\right) &=&e_{1}+e_{2}, \\
\phi \left( f_{2}\right) &=&-e_{1}-e_{2}+f_{2}\text{.}
\end{eqnarray*}

Direct verification reveals that $\phi ^{\ast }\left( B_{M}\oplus
B_{M}\right) =B_{N}\oplus B_{N}$, which is what we wanted to show.

Second, we show that 
\begin{equation*}
3\cdot \lbrack 1]=\left[ -1\right] +%
\begin{bmatrix}
2 & 1 \\ 
1 & 2%
\end{bmatrix}%
\text{.}
\end{equation*}

Let $V_{i}=\left( 
\mathbb{Z}
/4%
\mathbb{Z}
\cdot e_{i},B_{\left[ 1\right] }\right) $, $i=1,2,3$ be three copies of the
symmetric space associated to the matrix $\left( 1\right) $. Let $\phi
:V_{1}\oplus V_{2}\oplus V_{3}\rightarrow V_{1}\oplus V_{2}\oplus V_{3}$ be
the isomorphism given by 
\begin{eqnarray*}
\phi \left( e_{1}\right) &=&e_{1}+e_{2}+e_{3}, \\
\phi \left( e_{2}\right) &=&e_{1}+2e_{2}+e_{3}, \\
\phi \left( e_{3}\right) &=&e_{1}+e_{2}+2e_{3}\text{.}
\end{eqnarray*}

Direct verification reveals that $\phi ^{\ast }\left( B_{\left[ 1\right]
}\oplus B_{\left[ 1\right] }\oplus B_{\left[ 1\right] }\right) =B_{\left[ -1%
\right] }\oplus B_{M}$, which is what we wanted to show$.$

In the same fashion one shows that 
\begin{equation*}
3\cdot \lbrack -1]=\left[ 1\right] +%
\begin{bmatrix}
2 & 1 \\ 
1 & 2%
\end{bmatrix}%
.
\end{equation*}

This concludes the proof of the proposition.

\subsection{Proof of Proposition \protect \ref{structure3_prop}}

Let $\left[ V,B\right] \in \mathrm{GW}$. Using Proposition \ref%
{structure1_prop}, we can write 
\begin{equation*}
\left[ V,B\right] =n_{1}\cdot \left[ 1\right] +n_{2}\cdot \left[ -1\right]
+n_{3}\cdot 
\begin{bmatrix}
0 & 1 \\ 
1 & 0%
\end{bmatrix}%
+n_{4}\cdot 
\begin{bmatrix}
2 & 1 \\ 
1 & 2%
\end{bmatrix}%
.
\end{equation*}

Therefore 
\begin{eqnarray*}
4\cdot \left[ V,B\right] &=&4n_{1}\cdot \left[ 1\right] +4n_{2}\cdot \left[
-1\right] +4n_{3}\cdot 
\begin{bmatrix}
0 & 1 \\ 
1 & 0%
\end{bmatrix}%
+4n_{4}\cdot 
\begin{bmatrix}
2 & 1 \\ 
1 & 2%
\end{bmatrix}
\\
&=&n_{1}\cdot 
\begin{bmatrix}
2 & 1 \\ 
1 & 2%
\end{bmatrix}%
+n_{2}\cdot 
\begin{bmatrix}
2 & 1 \\ 
1 & 2%
\end{bmatrix}%
+4n_{3}\cdot \mathbf{0}+2n_{4}\cdot \mathbf{0} \\
&=&\left( n_{1}+n_{2}\right) 
\begin{bmatrix}
2 & 1 \\ 
1 & 2%
\end{bmatrix}%
=rk\left( V\right) \cdot \mathbf{4}\text{.}
\end{eqnarray*}

We have to explain the second and the forth equalities.

\textbf{The second equality. }In the second equality we used the following.
First, in $\mathrm{GW}$ 
\begin{eqnarray*}
\begin{bmatrix}
0 & 1 \\ 
1 & 0%
\end{bmatrix}
&=&\mathbf{0}\text{, } \\
\begin{bmatrix}
2 & 1 \\ 
1 & 2%
\end{bmatrix}%
+%
\begin{bmatrix}
2 & 1 \\ 
1 & 2%
\end{bmatrix}
&=&\mathbf{0,}
\end{eqnarray*}%
where in the second equality we used Proposition \ref{structure2_prop},
relation (\ref{rel1_eq}). Second, in $\mathrm{W}$ 
\begin{eqnarray*}
4\cdot \left[ 1\right] &=&\left[ 1\right] +3\cdot \left[ 1\right] =\left[ 1%
\right] +\left[ -1\right] +%
\begin{bmatrix}
2 & 1 \\ 
1 & 2%
\end{bmatrix}%
, \\
4\cdot \left[ -1\right] &=&\left[ -1\right] +3\cdot \left[ -1\right] =\left[
-1\right] +\left[ 1\right] +%
\begin{bmatrix}
2 & 1 \\ 
1 & 2%
\end{bmatrix}%
,
\end{eqnarray*}%
where, here, we used Proposition \ref{structure2_prop}, relations (\ref%
{rel2_eq}) and (\ref{rel3_eq}) respectively. Hence in $\mathrm{GW}$%
\begin{eqnarray*}
4\cdot \left[ 1\right] &=&%
\begin{bmatrix}
2 & 1 \\ 
1 & 2%
\end{bmatrix}%
, \\
4\cdot \left[ -1\right] &=&%
\begin{bmatrix}
2 & 1 \\ 
1 & 2%
\end{bmatrix}%
.
\end{eqnarray*}

\textbf{The Forth equality. }First, we note that the unique element in $%
\mathrm{GW}$ of order 2, which we denoted by $\mathbf{4}$ can be represented 
\begin{equation*}
\mathbf{4=}%
\begin{bmatrix}
2 & 1 \\ 
1 & 2%
\end{bmatrix}%
.
\end{equation*}

Second, since $rk\left( V\right) =n_{1}+n_{2}+2n_{3}+2n_{4}$ we have that $%
rk\left( V\right) =n_{1}+n_{2}$ $\left( \func{mod}2\right) $.

This concludes the proof of the proposition.

\subsection{Proof of Theorem \protect \ref{purity_thm}}

Since $G:\mathrm{W}\rightarrow 
\mathbb{C}
^{\times }$ is a morphism of monoids, it is enough to prove the assertion
for the generators $\left[ 1\right] ,\left[ -1\right] ,%
\begin{bmatrix}
0 & 1 \\ 
1 & 0%
\end{bmatrix}%
$ and $%
\begin{bmatrix}
2 & 1 \\ 
1 & 2%
\end{bmatrix}%
$, that is, we need to show%
\begin{eqnarray*}
\left \vert G\left( \left[ 1\right] \right) \right \vert &=&\sqrt{2}, \\
\left \vert G\left( \left[ -1\right] \right) \right \vert &=&\sqrt{2}, \\
\left \vert G\left( 
\begin{bmatrix}
0 & 1 \\ 
1 & 0%
\end{bmatrix}%
\right) \right \vert &=&2, \\
\left \vert G\left( 
\begin{bmatrix}
2 & 1 \\ 
1 & 2%
\end{bmatrix}%
\right) \right \vert &=&2.
\end{eqnarray*}

We calculate

\begin{equation*}
G\left( \left[ 1\right] \right) =\sum \limits_{x\in \mathbb{F}_{2}}e^{\frac{%
2\pi i}{4}x^{2}}=1+i=e^{\frac{2\pi i}{8}}\sqrt{2}\text{,}
\end{equation*}%
therefore, we get that $\left \vert G\left( \left[ 1\right] \right)
\right
\vert =\sqrt{2}$. Since $G\left( \left[ -1\right] \right) =\overline{%
G\left( \left[ 1\right] \right) }$, the assertion follows for $G\left( \left[
-1\right] \right) $ as well.

Now, we calculate 
\begin{equation*}
G\left( 
\begin{bmatrix}
0 & 1 \\ 
1 & 0%
\end{bmatrix}%
\right) =\sum \limits_{x,y\in \mathbb{F}_{2}}e^{\frac{2\pi i}{2}xy}=2.
\end{equation*}

Using the above equality and Proposition \ref{structure1_prop}, we can write 
\begin{equation*}
G\left( 
\begin{bmatrix}
2 & 1 \\ 
1 & 2%
\end{bmatrix}%
\right) ^{2}=G\left( 
\begin{bmatrix}
0 & 1 \\ 
1 & 0%
\end{bmatrix}%
\right) ^{2}=4\text{,}
\end{equation*}%
which implies that $\left \vert G\left( 
\begin{bmatrix}
2 & 1 \\ 
1 & 2%
\end{bmatrix}%
\right) \right \vert =2$.

This concludes the proof of the theorem.

\subsection{Proof of Theorem \protect \ref{GWitt_thm}}

Denote 
\begin{eqnarray*}
A &=&\left[ 1\right] +\left[ -1\right] , \\
B &=&%
\begin{bmatrix}
0 & 1 \\ 
1 & 0%
\end{bmatrix}%
\text{.}
\end{eqnarray*}

First, we show that $A,B\in I$. For this it is enough to show that 
\begin{eqnarray*}
G\left( A\right) &=&2, \\
G\left( B\right) &=&2.
\end{eqnarray*}

The first equation is clear since $G\left( \left[ -1\right] \right) =%
\overline{G\left( \left[ 1\right] \right) }$. The second equation is
obtained by direct calculation%
\begin{equation*}
G\left( B\right) =\sum \limits_{x,y\in \mathbb{F}_{2}}e^{\frac{2\pi i}{2}%
xy}=2\text{.}
\end{equation*}

Let us denote by $I^{\prime }\subset \mathrm{W}$ the submonoid generated by $%
A,B$.

Second, we show that $\mathrm{W}/I^{\prime }$ is a group, isomorphic to $%
\mathbb{Z}
/8%
\mathbb{Z}
$. First, note that Propositions \ref{structure1_prop} and \ref%
{structure2_prop} imply that the element $\left[ 1\right] $ generates the
group $\mathrm{W}/I^{\prime }$ and, in addition, $8\cdot \left[ 1\right] =%
\mathbf{0}$ in $\mathrm{W}/I^{\prime }$. This implies that the morphism of
monoids $%
\mathbb{Z}
/8%
\mathbb{Z}
\rightarrow \mathrm{GW}$, sending $1$ to $\left[ 1\right] $ is a surjection,
which in particular implies that $\mathrm{W}/I^{\prime }$ is a group - \ a
quotient group of $%
\mathbb{Z}
/8%
\mathbb{Z}
$.

In order to show that $\mathrm{W}/I^{\prime }$ is isomorphic to $%
\mathbb{Z}
/8%
\mathbb{Z}
$ it will be enough to show that $G\left( \left[ 1\right] \right) =\mu \sqrt{%
2}$, where $\mu $ is a primitive element in $\mu _{8}$. We calculate 
\begin{equation}
G\left( \left[ 1\right] \right) =\sum \limits_{x\in \mathbb{F}_{2}}e^{\frac{%
2\pi i}{4}x^{2}}=1+i=e^{\frac{2\pi i}{8}}\sqrt{2}\text{.}
\label{primitive_eq}
\end{equation}

Finally, we show that the canonical morphism $\mathrm{W}/I^{\prime
}\rightarrow \mathrm{GW}=\mathrm{W}/I$ is an isomorphism. This follows from
the fact that $G:\mathrm{W}/I^{\prime }\rightarrow 
\mathbb{C}
^{\times }/2^{%
\mathbb{Z}
}$ is an injection, which, in turns, follows from \ref{primitive_eq}.

This concludes the proof of the theorem.

\subsection{Proof of Proposition \protect \ref{gerbe-split_prop}}

We need to show that the category $\mathcal{W}^{s}$ is a gerbe with band $%
\mu _{2}$.

Clearly $\mathcal{W}^{s}$ is a groupoid. We need to show that every two
objects in $\mathcal{W}^{s}$ are isomorphic.

Consider two objects $\mathfrak{E}_{1}^{s},\mathfrak{E}_{2}^{s}\in \mathcal{W%
}^{s}$, \ that is 
\begin{eqnarray*}
\mathfrak{E}_{1}^{s} &=&(\widetilde{\mathfrak{E}}_{1},\mathfrak{E}%
_{1},\alpha _{1}), \\
\mathfrak{E}_{2}^{s} &=&(\widetilde{\mathfrak{E}}_{2},\mathfrak{E}%
_{2},\alpha _{2}),
\end{eqnarray*}%
where $\widetilde{\mathfrak{E}}_{i}\in \widetilde{\mathcal{W}}^{s}$, $%
\mathfrak{E}_{i}\in \mathcal{W}$ and $\alpha _{i}\in \mathrm{Mor}_{%
\widetilde{\mathcal{W}}}(\widetilde{S}(\widetilde{\mathfrak{E}}_{i}),\pi
^{\ast }(\mathfrak{E}_{i}))$, for $i=1,2$.

Since $\widetilde{\mathcal{W}}^{s}$ is a gerbe, there exists an morphism $%
\widetilde{f}\in \mathrm{Mor}_{\widetilde{\mathcal{W}}^{s}}(\widetilde{%
\mathfrak{E}}_{1},\widetilde{\mathfrak{E}}_{2})$. Since $\pi ^{\ast }$ is
full, there exists a morphism $f\in \mathrm{Mor}_{\mathcal{W}}(\mathfrak{E}%
_{1},\mathfrak{E}_{2})$ such that $\pi ^{\ast }\left( f\right) =\alpha
_{2}\circ \widetilde{S}(\widetilde{f})\circ \alpha _{1}^{-1}$. The pair $(%
\widetilde{f},f)$ is a morphism in $\mathrm{Mor}_{\mathcal{W}^{s}}(\mathfrak{%
E}_{1}^{s},\mathfrak{E}_{2}^{s})$.

We are left to show that $\mathrm{Mor}_{\mathcal{W}^{s}}(\mathfrak{E}^{s},%
\mathfrak{E}^{s})\simeq \mu _{2}$, for every $\mathfrak{E}^{s}\in \mathcal{W}%
^{s}$. Write $\mathfrak{E}^{s}=(\widetilde{\mathfrak{E}},\mathfrak{E},\alpha
)$. A morphism $(\widetilde{f},f)\in \mathrm{Mor}_{\mathcal{W}^{s}}(%
\mathfrak{E}^{s},\mathfrak{E}^{s})$ satisfies $\pi ^{\ast }\left( f\right)
=\alpha \circ \widetilde{S}(\widetilde{f})\circ \alpha ^{-1}$, therefore, $%
\pi ^{\ast }\left( f\right) $ is determined by $\widetilde{S}(\widetilde{f})$%
. Since $\widetilde{S}$ is faithful, $\pi ^{\ast }\left( f\right) $ is, in
fact, determined by $\widetilde{f}$. Finally, since $\pi ^{\ast }$ is
faithful, we obtain that $f$ is determined by $\widetilde{f}$. Therefore, we
get that 
\begin{equation*}
\mathrm{Mor}_{\mathcal{W}^{s}}(\mathfrak{E}^{s},\mathfrak{E}^{s})\simeq 
\mathrm{Mor}_{\widetilde{\mathcal{W}}^{s}}(\widetilde{\mathfrak{E}},%
\widetilde{\mathfrak{E}})\simeq \mu _{2}\text{.}
\end{equation*}

This concludes the proof of the proposition.

\subsection{Proof of Proposition \protect \ref{explicit1_prop}}

The proof proceeds along similar lines as the proof of Proposition \ref%
{explicit_prop}.

Let $\delta \in \mathcal{H}_{L}$ be the unique function on $H\left( V\right) 
$, supported on $Z\cdot L\subset H\left( V\right) $ and normalized such that 
$\delta \left( 0\right) =1$.

On the one hand, explicit computation reveals that $F_{\widetilde{N},%
\widetilde{L}}\left[ \delta \right] \left( 0\right) =1$. This implies that 
\begin{equation*}
C(\widetilde{N},\widetilde{M},\widetilde{L})=F_{\widetilde{N},\widetilde{M}%
}\circ F_{\widetilde{M},\widetilde{L}}\left[ \delta \right] \left( 0\right) 
\text{.}
\end{equation*}

On the other hand%
\begin{equation*}
F_{\widetilde{N},\widetilde{M}}\circ F_{\widetilde{M},\widetilde{L}}\left[
\delta \right] \left( 0\right) =\sum \limits_{m\in M}\sum \limits_{n\in
N}\delta \left( \tau _{M}\left( m\right) \cdot \tau _{N}\left( n\right)
\right) \text{,}
\end{equation*}%
where $\tau _{M}:M\rightarrow H\left( V\right) $ and $\tau _{N}:N\rightarrow
H\left( V\right) $ are the enhanced Lagrangians $\pi (\widetilde{M})$ and $%
\pi (\widetilde{L})$ respectively. The multiplication rule in $H\left(
V\right) $ implies that for every $m\in M$ and $n\in N$ we have that 
\begin{equation*}
\tau _{M}\left( m\right) \cdot \tau _{N}\left( n\right) =(m+n,\alpha _{%
\widetilde{M}}\left( m\right) +\alpha _{\widetilde{N}}\left( n\right) +\beta
\left( m,n\right) ).
\end{equation*}

Since $\delta \left( \tau _{M}\left( m\right) \cdot \tau _{N}\left( n\right)
\right) =0$ unless $m+n\in L$ we obtain that for every $m\in M$, the only
non-zero contribution to the sum $\sum \limits_{n\in N}\delta \left( \tau
_{M}\left( m\right) \cdot \tau _{N}\left( n\right) \right) $ comes from $%
n=-r^{L}\left( m\right) $.

Therefore, we obtain that%
\begin{eqnarray*}
F_{\widetilde{N},\widetilde{M}}\circ F_{\widetilde{M},\widetilde{L}}\left[
\delta \right] \left( 0\right) &=&\sum \limits_{m\in M}\delta \left( \tau
_{M}\left( m\right) \cdot \tau _{N}\left( -r^{L}\left( m\right) \right)
\right) \\
&=&\sum \limits_{m\in M}\delta (m-r^{L}\left( m\right) ,\alpha _{\widetilde{M%
}}(m)+\alpha _{\widetilde{N}}(-r^{L}(m))-\beta (m,r^{L}\left( m\right) ) \\
&=&\sum \limits_{m\in M}\delta (m-r^{L}\left( m\right) ,\alpha _{\widetilde{L%
}}(m-r^{L}\left( m\right) )+\alpha \left( m\right) ) \\
&=&\sum \limits_{m\in M}\psi \left( \alpha \left( m\right) \right) .
\end{eqnarray*}%
where $\alpha \left( m\right) =\alpha _{\widetilde{M}}(m)+\alpha _{%
\widetilde{N}}(r^{L}(m))-\beta (m,r^{L}\left( m\right) )-\alpha _{\widetilde{%
L}}(m-r^{L}\left( m\right) )$.

Now, we have 
\begin{eqnarray*}
\alpha _{\widetilde{M}}(m) &=&\widetilde{\beta }\left( \widetilde{m},%
\widetilde{m}\right) , \\
\alpha _{\widetilde{N}}(r^{L}(m)) &=&\widetilde{\beta }(r^{\widetilde{L}}(%
\widetilde{m}),r^{\widetilde{L}}(\widetilde{m})), \\
\alpha _{\widetilde{L}}(m-r^{L}\left( m\right) ) &=&\widetilde{\beta }(%
\widetilde{m}-r^{\widetilde{L}}(\widetilde{m}),\widetilde{m}-r^{\widetilde{L}%
}(\widetilde{m})), \\
\beta (m,r^{L}\left( m\right) ) &=&2\widetilde{\beta }(\widetilde{m},r^{%
\widetilde{L}}(\widetilde{m})),
\end{eqnarray*}%
where $\widetilde{m}$ is any element in $\widetilde{M}$ lying over $m\in M$.

Using the above, a direct calculation reveals that 
\begin{equation*}
\alpha (m)=\widetilde{\omega }(r^{\widetilde{L}}(\widetilde{m}),\widetilde{m}%
)\text{.}
\end{equation*}

Therefore, we obtain that \ 
\begin{eqnarray*}
C(\widetilde{N},\widetilde{M},\widetilde{L}) &=&\sum \limits_{m\in M}\psi (%
\widetilde{\omega }(r^{\widetilde{L}}(\widetilde{m}),\widetilde{m})) \\
&=&\sum \limits_{m\in M}\psi (\widetilde{\omega }_{\widetilde{L}}(\widetilde{%
m},\widetilde{m}))=G([\widetilde{M},tr(\widetilde{\omega }_{\widetilde{L}})])%
\text{.}
\end{eqnarray*}

This concludes the proof of the proposition.

\subsection{Proof of Proposition \protect \ref{disc_prop}}

Let us first consider the element $X^{\prime }\in \mathrm{W}\left( R\right) $
\begin{equation*}
X^{\prime }=[\widetilde{M},\widetilde{\omega }_{\widetilde{L}}]+[B_{%
\widetilde{N},\widetilde{M}}]+[B_{\widetilde{M},\widetilde{L}}]+[-B_{%
\widetilde{M},\widetilde{L}}],
\end{equation*}%
where we use the abbreviated notation $\left[ B\right] $ for $\left[ R^{n},B%
\right] $.

\begin{lemma}
\label{disc1_lemma}The discriminant $d\left( X^{\prime }\right) =$ $1$.
\end{lemma}

The proposition now follows from

\begin{lemma}
\label{disc2_lemma}Let $\left[ V,B\right] \in \mathrm{W}\left( R\right) $ be
such that $d\left( \left[ V,B\right] \right) =1$ then $d\left( \left[
V,tr\left( B\right) \right] \right) =1.$
\end{lemma}

This concludes the proof of the proposition.

\subsubsection{Proof of Lemma \protect \ref{disc1_lemma}}

Write 
\begin{eqnarray}
d\left( X^{\prime }\right) &=&d([\widetilde{M},\widetilde{\omega }_{%
\widetilde{L}}])\cdot d([B_{\widetilde{N},\widetilde{M}}])\cdot d([B_{%
\widetilde{M},\widetilde{L}}])\cdot d([-B_{\widetilde{N},\widetilde{L}}]) 
\notag \\
&=&\left( -1\right) ^{n}d([\widetilde{M},\widetilde{\omega }_{\widetilde{L}%
}])\cdot \widetilde{\omega }_{\wedge }(o_{\widetilde{M}},o_{\widetilde{N}%
})\cdot \widetilde{\omega }_{\wedge }(o_{\widetilde{L}},o_{\widetilde{M}%
})\cdot \widetilde{\omega }_{\wedge }(o_{\widetilde{L}},o_{\widetilde{N}}),
\label{disc1-1_eq}
\end{eqnarray}%
where, by construction, $d([B_{\widetilde{M},\widetilde{L}}])=\widetilde{%
\omega }_{\wedge }(o_{\widetilde{L}},o_{\widetilde{M}})$ for every $(%
\widetilde{M},\widetilde{L})\in \widetilde{U}_{2}$. We proceed to compute $d(%
\widetilde{[M},\widetilde{\omega }_{\widetilde{L}}])$.

First, we note 
\begin{equation*}
d([\widetilde{M},\widetilde{\omega }_{\widetilde{L}}])=\widetilde{\omega }%
_{\wedge }(r_{\wedge }^{\widetilde{L}}(o_{\widetilde{M}}),o_{\widetilde{M}}),
\end{equation*}%
where $r_{\wedge }^{\widetilde{L}}:\wedge ^{n}\widetilde{M}\rightarrow
\wedge ^{n}\widetilde{N}$ is the map induced from $r^{\widetilde{L}}:%
\widetilde{M}\rightarrow \widetilde{N}$.

\begin{lemma}
\label{disc3_lemma} We have 
\begin{equation}
\widetilde{\omega }_{\wedge }(r_{\wedge }^{\widetilde{L}}(o_{\widetilde{M}%
}),o_{\widetilde{M}})=\left( -1\right) ^{n}\frac{\widetilde{\omega }_{\wedge
}(o_{\widetilde{L}},o_{\widetilde{M}})\cdot \widetilde{\omega }_{\wedge }(o_{%
\widetilde{M}},o_{\widetilde{N}})}{\widetilde{\omega }_{\wedge }(o_{%
\widetilde{L}},o_{\widetilde{N}})}\text{.}  \label{disc1-2_eq}
\end{equation}
\end{lemma}

Substituting (\ref{disc1-2_eq}) in (\ref{disc1-1_eq}) we get 
\begin{equation*}
d\left( X^{\prime }\right) =\widetilde{\omega }_{\wedge }(o_{\widetilde{L}%
},o_{\widetilde{M}})^{2}\cdot \widetilde{\omega }_{\wedge }(o_{\widetilde{M}%
},o_{\widetilde{N}})^{2}\text{.}
\end{equation*}

This concludes the proof of the lemma.

\paragraph{Proof of Lemma \protect \ref{disc3_lemma}}

We know that $r_{\wedge }^{\widetilde{L}}(o_{\widetilde{M}})\in \wedge ^{n}%
\widetilde{N}$ is proportional to $o_{\widetilde{N}}$ since $\wedge ^{n}%
\widetilde{N}=R\cdot o_{\widetilde{N}}$, so we can write $r_{\wedge }^{%
\widetilde{L}}(o_{\widetilde{M}})=$ $a\cdot o_{\widetilde{N}}$, for some $%
a\in R^{\times }$.

On the one hand, since $r^{\widetilde{L}}\left( \widetilde{m}\right) -%
\widetilde{m}\in \widetilde{L}$, we have 
\begin{equation*}
\widetilde{\omega }_{\wedge }(r_{\wedge }^{\widetilde{L}}(o_{\widetilde{M}%
}),o_{\widetilde{L}})=\widetilde{\omega }_{\wedge }(o_{\widetilde{M}},o_{%
\widetilde{L}}).
\end{equation*}

On the other hand, since $r_{\wedge }^{\widetilde{L}}(o_{\widetilde{M}})=$ $%
a\cdot o_{\widetilde{N}}$, we have 
\begin{equation*}
\widetilde{\omega }_{\wedge }(r_{\wedge }^{\widetilde{L}}(o_{\widetilde{M}%
}),o_{\widetilde{L}})=a\cdot \widetilde{\omega }_{\wedge }(o_{\widetilde{N}%
},o_{\widetilde{L}}).
\end{equation*}

Combining the above two equations we get 
\begin{equation*}
a=\frac{\widetilde{\omega }_{\wedge }(o_{\widetilde{M}},o_{\widetilde{L}})}{%
\widetilde{\omega }_{\wedge }(o_{\widetilde{N}},o_{\widetilde{L}})}=\frac{%
\widetilde{\omega }_{\wedge }(o_{\widetilde{L}},o_{\widetilde{M}})}{%
\widetilde{\omega }_{\wedge }(o_{\widetilde{L}},o_{\widetilde{N}})}.
\end{equation*}

Finally, write 
\begin{eqnarray*}
\widetilde{\omega }_{\wedge }(r_{\wedge }^{\widetilde{L}}(o_{\widetilde{M}%
}),o_{\widetilde{M}}) &=&a\cdot \widetilde{\omega }_{\wedge }(o_{\widetilde{N%
}},o_{\widetilde{M}}) \\
&=&\frac{\widetilde{\omega }_{\wedge }(o_{\widetilde{L}},o_{\widetilde{L}})}{%
\widetilde{\omega }_{\wedge }(o_{\widetilde{L}},o_{\widetilde{N}})}\cdot 
\widetilde{\omega }_{\wedge }(o_{\widetilde{N}},o_{\widetilde{M}}) \\
&=&\left( -1\right) ^{n}\frac{\widetilde{\omega }_{\wedge }(o_{\widetilde{L}%
},o_{\widetilde{M}})\cdot \widetilde{\omega }_{\wedge }(o_{\widetilde{M}},o_{%
\widetilde{N}})}{\widetilde{\omega }_{\wedge }(o_{\widetilde{L}},o_{%
\widetilde{N}})}\text{.}
\end{eqnarray*}

This concludes the proof of the lemma.

\subsubsection{Proof of Lemma \protect \ref{disc2_lemma}}

It is easy to verify that the discriminants $d\left( \left[ V,tr\left(
B\right) \right] \right) $ and $d\left( \left[ V,B\right] \right) $ are
related as follows%
\begin{equation*}
d\left( \left[ V,tr\left( B\right) \right] \right) =N\left( d\left( \left[
V,B\right] \right) \right) \cdot d\left( \left[ R,tr\right] \right)
^{rk_{R}\left( V\right) }\text{.}
\end{equation*}%
where $N:R^{\times }\rightarrow \left( 
\mathbb{Z}
/4%
\mathbb{Z}
\right) ^{\times }$ is the norm map. Since we assume that $d\left( \left[ V,B%
\right] \right) $ is a square in $R^{\times }$ this implies that $N\left(
d\left( \left[ V,B\right] \right) \right) $ is a square in $\left( 
\mathbb{Z}
/4%
\mathbb{Z}
\right) ^{\times }$ therefore it is enough to prove that 
\begin{equation*}
d\left( \left[ R,tr\right] \right) =1,
\end{equation*}%
where $tr:R\times R\rightarrow 
\mathbb{Z}
/4%
\mathbb{Z}
$ denote the trace form $tr\left( x,y\right) =tr\left( x\cdot y\right) $.

We prove a more general assertion. Let $K/F$ be an unramified extension of
local fields lying over $%
\mathbb{Q}
_{2}$. Denote $R_{K}=\mathcal{O}_{K}/\mathfrak{m}_{K}^{2}$ and $R_{F}=%
\mathcal{O}_{F}/\mathfrak{m}_{F}^{2}$.$\ $We claim that 
\begin{equation*}
d\left( \left[ R_{K},tr_{K/F}\right] \right) =1\text{.}
\end{equation*}%
where $tr_{K/F}:R_{K}\times R_{K}\rightarrow R_{F}$ denote the relative
trace map. First, we prove the above in two particular cases. The first case
is when $\left[ K:F\right] $ is odd; the second case is when $\left[ K:F%
\right] =2^{l}$, for some $l\in 
\mathbb{N}
$.

\textbf{Case 1. }Assume $d=\left[ K:F\right] $ is odd. In this case we claim
that already $d\left( \left[ \mathcal{O}_{K},tr_{K/F}\right] \right) =1$ in $%
\mathcal{O}_{F}^{\times }/\mathcal{O}_{F}^{\times 2}$.

We can present $K$ in the form $K=F\left( \alpha \right) $ such that $%
1,\alpha ,\alpha ^{2},...,\alpha ^{d-1}$ is a basis of $\mathcal{O}_{K}$
over $\mathcal{O}_{F}$. Explicit computation reveals that 
\begin{eqnarray}
d\left( \left[ \mathcal{O}_{K},tr_{K/F}\right] \right) &=&\Delta \left(
1,\alpha ,\alpha ^{2},...,\alpha ^{d-1}\right)  \notag \\
&=&\tprod \limits_{0\leq i<j\leq d-1}\left( Fr^{j}\alpha -Fr^{i}\alpha
\right) ^{2},  \label{disc2-1_eq}
\end{eqnarray}%
where $Fr\in G=Gal\left( K/F\right) $ is the Frobenius automorphism. Denote 
\begin{equation*}
D=\tprod \limits_{0\leq i<j\leq d-1}\left( Fr^{j}\alpha -Fr^{i}\alpha
\right) .
\end{equation*}

We claim that $D\in \mathcal{O}_{F}^{\times }$. Clearly, $D\in \mathcal{O}%
_{K}^{\times }$. In addition, for every $g\in G$, we have that $gD=\sigma
\left( g\right) D$, where $\sigma :\Sigma _{d}\rightarrow \left \{ \pm
1\right \} $ is the sign homomorphism of the permutation group and we use
the injection $G\hookrightarrow \Sigma _{d}$. Since $K/F$ is an unramified
extension this implies that $G$ is cyclic and since by assumption $%
\left
\vert G\right \vert $ is odd, it must be that $\sigma _{|G}=1$.
Hence, $D $ is invariant under the action of $G$ therefore $D\in \mathcal{O}%
_{F}^{\times }$.

Combining this with (\ref{disc2-1_eq}) we get that $d\left( \left[ \mathcal{O%
}_{K},tr_{K/F}\right] \right) =1$. This concludes the proof of the assertion
in this case.

\textbf{Case 2.} Assume $d=\left[ K:F\right] =2^{l}$, for some $l\in 
\mathbb{N}
$. Consider a series of intermediate extensions interpolating between $K$
and $F$ 
\begin{equation*}
K=K_{l}-K_{l-1}-...-K_{1}-K_{0}=F,
\end{equation*}%
where $\left[ K_{i}:K_{i-1}\right] =2$, for every $i=1,..,l$.

Denote $R_{i}=\mathcal{O}_{K_{i}}/\mathfrak{m}_{K_{i}}^{2}$, let $%
tr_{i/i-1}:R_{i}\rightarrow R_{i-1}$ and $tr_{i}:R_{i}\rightarrow
R_{0}=R_{F} $ denote the corresponding trace maps and finally, let $%
N_{i}:R_{i}^{\times }\rightarrow R_{0}^{\times }=R_{F}^{\times }$ denote the
norm map.

It is easy to verify that the discriminants $d\left( \left[ R_{i},tr_{i}%
\right] \right) $ and $d\left( \left[ R_{i-1},tr_{i-1}\right] \right) $ are
related as follows 
\begin{equation*}
d\left( \left[ R_{i},tr_{i}\right] \right) =N_{i-1}\left( d\left( \left[
R_{i},tr_{i/i-1}\right] \right) \right) \cdot d\left( \left[ R_{i-1},tr_{i-1}%
\right] \right) \text{.}
\end{equation*}

It is enough to prove that $d\left( \left[ R_{i},tr_{i/i-1}\right] \right)
=1 $ in $R_{i-1}^{\times }/R_{i-1}^{\times 2}$.

We can assume that $K_{i}=K_{i-1}\left( \alpha \right) $ where $\alpha
^{2}+\alpha +a=0$, for some $a\in \mathcal{O}_{K_{i-1}}$. Now, we have 
\begin{equation*}
d\left( \left[ \mathcal{O}_{K_{i}},tr_{i/i-1}\right] \right) =\Delta \left(
1,\alpha \right) =\left( \alpha _{1}-\alpha _{2}\right) ^{2},
\end{equation*}%
where $\alpha _{1},\alpha _{2}$ are the two roots of the polynomial $%
x^{2}+x+a$. Explicit computation reveals that $\Delta \left( 1,\alpha
\right) =1-4a$ which implies that $d\left( \left[ R_{i},tr_{i/i-1}\right]
\right) =1$.

This concludes the proof of the assertion in this case.

The statement for a general extension $K/F$ now follows easily: Let $E$ be
an intermediate extension $K-E-F$ such that $\left[ K:E\right] $ is odd and $%
\left[ E:F\right] =2^{l}$ and use 
\begin{equation*}
d\left( \left[ R_{K},tr_{K/F}\right] \right) =N_{E/F}\left( d\left( \left[
R_{K},tr_{K/E}\right] \right) \right) \cdot d\left( \left[ R_{F},tr_{E/F}%
\right] \right) .
\end{equation*}

This concludes the proof of the lemma.

\subsection{Proof of Proposition \protect \ref{vanishing_prop}}

Write $X\in \mathrm{W}$ in the form 
\begin{equation*}
X=n_{1}\left[ 1\right] +n_{2}\left[ -1\right] +n_{3}%
\begin{bmatrix}
0 & 1 \\ 
1 & 0%
\end{bmatrix}%
+n_{4}%
\begin{bmatrix}
2 & 1 \\ 
1 & 2%
\end{bmatrix}%
\text{.}
\end{equation*}

We have $rk\left( X\right) =n_{1}+n_{2}+2n_{3}+2n_{4}$. On the one hand,
since $4|rk\left( X\right) $ we have 
\begin{equation}
n_{1}+n_{2}+2n_{3}+2n_{4}\equiv 0\text{ }\func{mod}4\text{.}
\label{vanish1_eq}
\end{equation}

On the other hand, $d\left( X\right) =\left( -1\right) ^{n_{1}+n_{2}+n_{3}}$%
, where we use here the facts that 
\begin{equation*}
d\left( 
\begin{bmatrix}
0 & 1 \\ 
1 & 0%
\end{bmatrix}%
\right) =d\left( 
\begin{bmatrix}
2 & 1 \\ 
1 & 2%
\end{bmatrix}%
\right) =-1.
\end{equation*}

Therefore, since $d\left( X\right) =1$ we have $n_{2}+n_{3}+n_{4}\equiv 0$ $%
\func{mod}2$, which implies that 
\begin{equation}
2n_{2}+2n_{3}+2n_{4}\equiv 0\func{mod}4.  \label{vanish2_eq}
\end{equation}

Substracting (\ref{vanish2_eq}) from (\ref{vanish1_eq}) we get that $n_{1}$ $%
\equiv n_{2}$ $\func{mod}4$ which implies that 
\begin{equation}
2n_{1}\equiv 2n_{2}\func{mod}8\text{.}  \label{vanish3_eq}
\end{equation}

Finally, in $\mathrm{GW}$, we can write 
\begin{eqnarray*}
2X &=&2n_{1}\left[ 1\right] +2n_{2}\left[ -1\right] +2n_{3}%
\begin{bmatrix}
0 & 1 \\ 
1 & 0%
\end{bmatrix}%
+2n_{4}%
\begin{bmatrix}
2 & 1 \\ 
1 & 2%
\end{bmatrix}
\\
&=&2n_{1}\left[ 1\right] +2n_{2}\left[ -1\right] =\mathbf{0,}
\end{eqnarray*}%
where the second equality follows from the facts that $2%
\begin{bmatrix}
2 & 1 \\ 
1 & 2%
\end{bmatrix}%
=2%
\begin{bmatrix}
0 & 1 \\ 
1 & 0%
\end{bmatrix}%
$ in $\mathrm{W}$ (Theorem \ref{structure1_prop}) and $%
\begin{bmatrix}
0 & 1 \\ 
1 & 0%
\end{bmatrix}%
=\mathbf{0}$ in $\mathrm{GW}$ (Theorem \ref{GWitt_thm}) and the third
equality follows from (\ref{vanish3_eq}), the relation $\left[ 1\right] +%
\left[ -1\right] =\mathbf{0}$ (Theorem \ref{GWitt_thm}) and the fact that
for every $X\in \mathrm{GW}$, $8X=\mathbf{0}$ (Theorem \ref{GWitt_thm}).

This concludes the proof of the proposition.


\begin{thebibliography}{9}
\bibitem{B} Brylinski J-L., Central extensions and reciprocity laws. \textit{%
Cahiers Topologie G\'{e}om. Diff\'{e}rentielle Cat\'{e}g. 38, }no. 3 (1997),
193-215.\textit{\ }

\bibitem{BL} Blasco L., Paires duales r\'{e}ductives en caract\'{e}ristique
2. \textit{M\'{e}moires de la S.M.F. 2}$^{e}$ \textit{s\'{e}rie, tome 52}
(1993), p. 1-73.

\bibitem{EG} Etingof P. and Gelaki S., Isocategorical groups. \textit{%
Internat. Math. Res. Notices, }no. 2 (2001)\textit{\ }59-76\textit{.}

\bibitem{G} G\'{e}rardin P., Weil representations associated to finite
fields. \textit{\ J. Algebra 46}, no. 1, (1977), 54-101.

\bibitem{H} Howe R., On the character of Weil's representation. \textit{%
Trans. Amer. Math. Soc, }177 (1973), 287--298.

\bibitem{W} Weil A., Sur certains groupes d'operateurs unitaires. \textit{%
Acta Math. }111 (1964), 143-211.
\end{thebibliography}
\end{document}